\documentclass[12pt]{article}

\usepackage{upgreek}
\usepackage{graphicx}
\usepackage{float}%

\usepackage{hyperref}
\hypersetup{colorlinks = true} 

\usepackage{amsmath}

\usepackage{amsfonts}
\usepackage{amssymb}
\usepackage{xypic}
\usepackage{enumerate}
\usepackage{subfigure}
\usepackage{amsthm}
\usepackage{mathrsfs}

\theoremstyle{plain}
\newtheorem{thm}{Theorem}[section]
\newtheorem{theorem}[thm]{Theorem}

\newtheorem{corollary}[thm]{Corollary}

\newtheorem{lemma}[thm]{Lemma}

\newtheorem{prop}[thm]{Proposition}
\newtheorem{definition}[thm]{Definition}

\renewcommand{\phi}{\varphi}

\newcommand{\Div}{\mathrm{Div}}

\renewcommand{\O}{\mathscr{O}}
\newcommand{\C}{\mathbb{C}}

\newcommand{\Q}{\mathbb{Q}}
\newcommand{\R}{\mathbb{R}}
\newcommand{\Z}{\mathbb{Z}}

\newcommand{\rank}{\mathrm{rank}}

\newcommand{\Pic}{\mathrm{Pic}}
\newcommand{\hPic}{\widehat{\mathrm{Pic}}}
\newcommand{\im}{\mathrm{im}}
\renewcommand{\div}{\mathrm{div}}

\newcommand{\hDiv}{\widehat{\mathrm{Div}}}
\newcommand{\ds}{\displaystyle}

\newcommand{\can}{\mathrm{can}}
\newcommand{\ad}{\mathrm{ad}}
\newcommand{\Ar}{\mathrm{Ar}}
\renewcommand{\mod}{\mathrm{mod}}
\newcommand{\mult}{\mathrm{mult}}

\title{Asymptotic Behavior of the Zhang--Kawazumi's $ \phi $-invariants}
\author{Yinchong Song}

\begin{document}
	
\maketitle

\tableofcontents

\section{Introduction}
In 2009, Shou-Wu Zhang \cite{Zha10} introduced the $ \phi $-invariant for a smooth projective curve $ C $ over a local field $ K $. If $ K $ is non-archimedean, the $ \phi $-invariant of $ C $ depends only on the reduction graph $ \Gamma(C) $ of $ C $, and we can define the $ \phi $-invariants for any polarized metrized graphs. If $ K $ is archimedean, it is defined by the Laplacian operator for the Arakelov metric. The archimedean case is independently introduced by Nariya Kawazumi in \cite{Kaw08,Kaw09}. 

In this paper, we study some asymptotic properties of the $ \phi $-invariants of different kinds. More precisely, this paper includes the following two parts.

(1) Prove the continuity of Zhang's $ \phi $-invariants for degenerating graphs. 

(2) Show that the continuity in (1) induces an adelic divisor $ \Phi $ on the moduli space of Riemann surfaces of a special kind, and then give an asymptotic expression of the Zhang--Kawazumi $ \phi $-invariants for Riemann surfaces near the boundary of the moduli space. 

A similar asymptotic behavior of the $ \phi $-invariants for Riemann surfaces was conjectured by Robin de Jong in \cite[Conj.1.2]{dJo15}. His conjecture for one-parameter family was proved by Wilms \cite{Wil} and by Robin de Jong himeself and Farbod Shokriehin \cite{JS}. Our results hold for high-dimensional base. 

To prove the first continuity theorem, we use the theory of metrized graphs, and show that many other invariants are also continuous for degenerating graphs. It is also used in the second part.

The second theorem is based on the recent work of adelic line bundles of Yuan--Zhang \cite{YZ} and the globalization of the $ \phi $-invariants by Yuan \cite{Yua}, which gives a limit version of divisors and line bundles, together with an intersection theory of them. We need to combine the theories of graphs and adelic line bundles together. 

\subsection{Zhang's $ \phi $-invariants for Graphs}
Let $ (\Gamma,q) $ be a polarized metrized graph. See \S \ref{PolarizedMetrizedGraph} for details. Its $ \phi $-invariant is defined as
\begin{align*}
\phi(\Gamma,q) &= -\frac{1}{4} \ell(\Gamma) + \frac{1}{4} \int_{\Gamma} g_{\mu_{\ad}}(x,x) ((10g+2)\mu_{\mathrm{ad}} -\delta_{K}),
\end{align*}  
where $ \ell(\Gamma) $ is the total length of $ \Gamma $, $ K $ is the canonical divisor on $ \Gamma $, $ \mu_{\ad} $ is admissible measure on $ \Gamma $, and $ g_{\mu_{\mathrm{ad}}} $ is the admissible Green function on $ \Gamma\times \Gamma $.

The $ \phi $-invariants for graphs were studied by Faber \cite{Fab} and Cinkir \cite{Cin} before. Faber showed that the $ \phi $-invariants is a rational function of the lengths of edges of the graph, and he proved the positivity of the $ \phi $-invariants of genus $ g=2,3,4 $. Cinkir showed the positivity of the $ \phi $-invariants for all genera, which can be used in proving the Bogomolov conjecture for function fields. 

These two papers both consider the change of the $ \phi $-invariants when varying graphs, especially when the length of some edges tend to zero. In \cite[Prop. 5.4]{Fab}, Faber showed the continuity of the $ \phi $-invariants while the underlying topological space is a homotopy. Further, they both use small circles to replace the polarization in \cite[Lemma 5.14]{Fab} and \cite[Prop.4.16]{Cin}. 

The continuity property of the $ \phi $-invariants (and some other invariants including $ \varepsilon $-invariants) was also studied in \cite[Lem. 1.3]{Yam}, \cite[\S 1.4]{Wan}.

In our paper, we will slightly generalize their continuity results. 

Let $ (\Gamma,q) $ be a polarized graph, and let $ e_1,\ldots,e_m $ be edges of $ \Gamma $. Then for each $ L_1,\ldots,L_m>0 $, there is a metrized graph $ (\Gamma,q;L_1,\ldots,L_m) $, such that the length of $ e_i $ is $ L_i $. Its $ \phi $-invariant is denoted by $ \phi(\Gamma,q;L_1,\ldots,L_m) $, which is viewed as a function on $ (L_i)_{1\leqslant i\leqslant m}\in \R_{>0}^m $, and is called the $ \phi $-function.  

The space $ \R_{>0}^m $ is a ``moduli space" of (polarized) metrized graphs with fixed underlying graph. Then $ \R_{\geqslant0}^m $ is the closure of $ \R_{>0}^m $ in the Euclidean space $ \R^m $. We will give a moduli interpretation of $ \R_{\geqslant 0}^m $, such that each point in $ \R_{\geqslant 0}^m $ parametrizes a (polarized) metrized graph. This is called the degeneration of (polarized) metrized graph. We give the precise definition in \S \ref{DegenerationofGraphs}. 

Now we consider the $ \phi $-function on $ \R_{>0}^m $. Faber showed that this function is a rational function on $ \R_{>0}^m $ in \cite[Prop. 4.6]{Fab}. Our main theorem is as follows.
\begin{theorem} \label{PhiforGraph}
	The $ \phi $-function $ \phi(\Gamma,q;L_1,\ldots,L_m) $ extends to a continuous function on $ \R_{\geqslant0}^m $, which is also denoted by $ \phi(\Gamma,q;L_1,\ldots,L_m) $. 
	
	Further, the value $ \phi(\Gamma,q;L_1,\ldots,L_m) $ on the boundary is just the $ \phi $-invariant of the polarized metrized graph parametrized by that point. 
\end{theorem}

The proof of this theorem uses a careful analysis on graph theory. First, we review some basic properties of metrized graphs in \S \ref{PolarizedMetrizedGraph}. Then we will define a notion of degeneration of polarized metrized graphs, and show that the $ \phi $-invariants is continuous under degeneration in \S\ref{ContinuityProperties}. We will show that many functions, including the voltage function $ j(x,y,z) $ on $ \Gamma^3 $ and the admissible green function $ g_{\mu_{\ad}}(x,y) $ on $ \Gamma^2 $, are continuous under degeneration, and then deduce our continuity of the $ \phi $-invariants. 

In the second part we will show that the continuity of the $ \phi $-invariant characterizes the adelic divisor $ \Phi $ introduced by \cite[\S 3.3.2]{Yua}, which is a globalization of the non-archimedean $ \phi $-invarianst; and the $ \phi $-function for graphs determines the behavior of the $ \phi $-function for Riemann surfaces.

\subsection{Zhang--Kawazumi's $ \phi $-invariants for Riemann Surfaces}
The $ \phi $-invariants for Riemann surfaces were independently introduced by Kawazumi \cite{Kaw08,Kaw09} and Zhang \cite{Zha10}. Let $ C $ be a compact Riemann surface of genus $ g\geqslant 2 $. Its $ \phi $-invariant is defined by
\begin{align*}
\phi(C) = -\int_{(C\times C)} g_{\Ar}(x,y) c_1(\bar\O(\Delta))^2,
\end{align*}
where $ g_{\Ar}(x,y) $ is the Arakelov Green's function \cite{Ara} on $ C\times C $, and $ \Delta $ is the diagonal divisor on $ C\times C $. The metric on $ \bar{\O}(\Delta) $ is given by 
\begin{align*}
-\log \|1\|(x,y)= g_{\Ar}(x,y)
\end{align*}
for all $ (x,y)\in (C\times C)\backslash \Delta $. Here $ 1 $ denotes the canonical section of $ \O(\Delta) $.

Let $ S $ be an irreducible and reduced scheme or an analytic space. Let $ \pi:X\to S $ be a projective, generically smooth family of stable curves of genus $ g\geqslant 2 $, and let $ U\subset S $ parametrize the smooth curves. Then if $ s\in U $, the fiber $ X_s $ is a Riemann surface, and we get an invariant $ \phi(s)=\phi(X_s) $. Now the $ \phi $-invariant becomes a function on $ U $, and we are interested in the behavior of the $ \phi $-function when $ s $ approaches $ S\backslash U $.

Any family of stable curves $ \pi:X\to S $ is analytic locally (or \'etale locally) pull-back of the Kuranishi family, see \cite[\S XI.4]{ACGH}. Let $ s_0\in S\backslash U $ and $ X_{0} $ be the fiber over $ s_0 $ with $ r $ nodes $ x_1,\ldots,x_r\in X_0 $. Let $ \pi':Y\to T $ be the Kuranishi family of $ X_{0} $. By shrinking $ S $ if necessary, there is a unique holomorphic map $ u:S\to T $ such that the family $ X\to S $ is the pull-back of $ Y\to T $ via $ u $. 

By shrinking $ T $ if necessary, we may assume that $ T $ is a polydisk of dimension $ n=3g-3 $, centered at $ 0 $, and with a local coordinate $ (t_1,\cdots,t_n) $ in such a way that for $ 1\leqslant i\leqslant r $, the locus $ t_i=0 $ parameterizes deformations which are locally trivial at the node $ x_i $. Then the singular locus on $ T $ is defined by the equation $ t_1\cdots t_r=0 $. By pull-back, let $ u_i=t_i\circ u $ be functions on $ S $, then the singular locus $ Z=S\backslash U $ on $ S $ is defined by the equation $ u_1\cdots u_r=0 $.

For $ s\in S $, let $ X_s $ be the corresponding curve over $ s $, and if $ s\in U $, let $ \phi(X_s) $ be the $ \phi $-invariant of the Riemann surface $ X_s $.

On the other hand, let $ (\Gamma,q) $ be the dual graph of $ X_{0} $, and let 
$$ \tilde{\phi}(\Gamma,q;L_1,\ldots,L_r) $$
be the $ \phi $-function of the polarized graph $ (\Gamma,q) $. 
\begin{theorem}\label{PhiForRiemannSurface1}
	Let $ \delta $ be a real number with $ 0<\delta<1  $. Then the $ \phi $-function on $ S $ has the following asymptotic behavior
	\begin{align*}
	\phi(X_s)= \tilde{\phi}(\Gamma,q; -\log |u_1(s)|,\ldots,-\log |u_r(s)|) + o\left(\sum_{i=1}^r\left(-\log |u_i(s)| \right) \right)
	\end{align*}
	for $ s\in U $, $ |u_i(s)|\leqslant \delta $, and  $ \min_{1\leqslant i\leqslant r} |u_i(s)|\to 0 $. 
\end{theorem}

The $ \phi $-invariants of Riemann surfaces have also been widely studied before. 

When $ \dim S=1 $, Robin de Jong \cite{dJo14} studied the asymptotic behavior of the $ \phi $-invariants when $ X_s $ degenerates to a stable curve with only one node. Later, in \cite{dJo15}, he showed a similar results for high-dimensional base when the genus $ g=2 $, and for one-dimensional base for hyperelliptic curves of any genus $ g\geqslant 2 $. He conjectured that the $ \phi $-invariants have a simple form in all cases, which is proved by Wilms \cite{Wil} and by de Jong himself and Shokrieh \cite{JS} for one-dimensional base recently. They use a careful analysis of Arakelov-Green's function and Faltings's delta function. All their results have a better error term $ O(1) $.

Theorem \ref{PhiForRiemannSurface1} confirms de Jong's conjecture for high-dimensional base with a different error term. Under additional assumptions, we show that the error term in Theorem \ref{PhiForRiemannSurface1} indeed can be refined by $ O(1) $.

Let $ X $ be a stable curve. A node $ p\in X $ is called a separating node of type $ (i,g-i) $ for $ 1\leqslant i\leqslant [g/2] $ if the partial normalization of $ X $ at $ p $ is the union of two connected component of genera $ i $ and $ g-i $. 
It is called a non-separating node of type $ 0 $ if the partial normalization is connected. 

If $ X_0 $ has no non-separating nodes, then the $ \phi $-function $ \tilde{\phi}(\Gamma,q;L_1,\ldots,L_r) $ has an explicit form, and we have a better error term as follows. 

\begin{theorem}\label{PhiForRiemannSurface2} 
	Use the same notations as in Theorem \ref{PhiForRiemannSurface1}. Assume that $ X_0 $ has no non-separating nodes. Let $ (g_i,g-g_i) $ be the type of the node $ x_i $. Then the $ \phi $-invariant $ \phi(X_s) $ has the following asymptotic behavior
	\begin{align*}
	\phi(X_s) = -\sum_{i=1}^r \frac{2g_i(g-g_i)}{g} \log |u_i(s)| +O(1),
	\end{align*}
	where $ s\in U $ near $ s_0 $, and $ O(1) $ mean a continuous function on $ U $ which extends continuously to $ S $.
\end{theorem}

These two theorems are proved separately. Both proofs are based on Yuan--Zhang's recent work of adelic line bundles in \cite{YZ} and Yuan's globalization of the $ \phi $-invariants in \cite{Yua}. The notion of adelic divisors and line bundles is a limit version of usual (hermitian) divisors and line bundles. Yuan showed that the $ \phi $-invariants for both non-archimedean and archimedean cases can be globalized to an adelic divisor $ \bar{\Phi}=(\Phi,g_{\bar{\Phi}}) $ on the moduli space of smooth curves in \cite[\S 3.3.2]{Yua}. Then he got the positivity result for the $ \phi $-invariants of Riemann surfaces, and proved the uniform Bogomolov conjecture. 

In our paper, we will follow Yuan's idea to study the adelic divisor $ \bar{\Phi} $ on the moduli space of smooth curves. We will review basic definitions and properties of adelic divisors in section \ref{AdelicDivisors} and Yuan's globalization of the $ \phi $-invariants $ \bar{\Phi}=(\Phi,g_{\bar{\Phi}}) $ in \S \ref{GlobalizationofPhi}. Then in \S \ref{Proof2}, we will show that this adelic divisor $ \bar{\Phi} $ extends to the moduli space of stable curves without non-separating nodes, and thus get Theorem \ref{PhiForRiemannSurface2}. 

Theorem \ref{PhiForRiemannSurface1} will be proved differently. Recall that Zhang \cite{Zha93} gave two equivalent definitions of the admissible line bundles, based on the graph theory or Tate's limiting argument of abelian varieties. Yuan \cite{Yua} used abelian varieties to generalize the admissible line bundles to family of smooth curves, and then he got the globalization of the $ \Phi $-invariants. In our paper, we will show that Theorem \ref{PhiforGraph} can be used to give another definition of the underlying divisor $ \Phi $.  In \S \ref{Special}, we define special adelic divisors to be those satisfying a similar property. Then we give a criterion for special adelic divisors, and shows that for these adelic divisors, there exists a Green function of particular form. Then we get Theorem \ref{PhiForRiemannSurface1} in \S \ref{Proof1}. The error term in Theorem \ref{PhiForRiemannSurface1} is just the error term of two Green functions of the same underlying divisor as in \cite[Thm. 3.6.4]{YZ}.  

In \S \ref{OtherResults}, we will show that many other adelic divisors induced by the admissible line bundles are special. In particular, we can get a multi-dimensional asymptotic of Arakelov Green's function and the canonical metrics, but the error terms are just like $ o(-\sum \log |t_i|) $. In \cite[Thm. 1.1]{BGHJ}, José Burgos Gil, David Holmes and Robin de Jong get a multi-dimensional asymptotic of the canonical admissible metric on the theta line bundle on a family of complex polarized abelian varieties, and their error term is much better than ours. They get $ O(1) $, which is further continuous away from the singularities of the boundary divisor.

\subsubsection*{Acknowledgment}
I would like to thank Xinyi Yuan for explaining to me the problem of the $ \phi $-invariants, and for his continued guidance and useful suggestion to this paper. I would thank Jiawei Yu for his helpful feedback on earlier version of this paper. I would also thank Robin de Jong for his valuable comments.

\section{Degeneration of Polarized Graphs and Continuity Properties}

In this section, we will define the degeneration of graphs, and show that many functions and invariants are continuous under degeneration.

\subsection{Preliminaries on Metrized Graphs}\label{PolarizedMetrizedGraph}
We first review the theory of polarized metrized graphs and then give the precise definition of Zhang's admissible metrics and $ \phi $-invariants. The theory of metrized graphs can be found in \cite{CR,BR,BF}, and Zhang's admissible metric and $ \phi $-invariants can be found in \cite{Zha93,Zha10}. 

\subsubsection*{Metrized Graphs}

A metrized graph $ \Gamma $ is a finite connected graph together with a distinguished parameterization of each edge. Once the graph is metrized, there is a Lebesgue measure on $ \Gamma $, denoted by $ dx $. 

For any $ p\in \Gamma $, the number of directions emanating from $ p $ is called the valence of $ p $, and is denoted by $ v(p) $. By definition, there are only finitely many $ p\in \Gamma $ with $ v(p)\neq 2 $.

The set of vertices in $ \Gamma $ is denoted by $ V(\Gamma) $. We require that $ V(\Gamma) $ is non-empty and contains all $ p\in \Gamma $ with $ v(p)\neq 2 $, such that $ \Gamma\backslash V(\Gamma) $ is disjoint unions of open intervals. For a given metrized graph $ \Gamma $, one can enlarge the vertex set $ V(\Gamma) $ by considering arbitrarily many valence $ 2 $ points as vertices. 

By an edge of $ \Gamma $, we mean a closed line segment with end points in $ V(\Gamma) $. The edge set of $ \Gamma $ is denoted by $ E(\Gamma) $. 

The genus of a graph $ \Gamma $ is defined as 
\begin{align*}
g(\Gamma)=\rank H^1(\Gamma,\Z).
\end{align*}
Since $ \Gamma $ is connected, $ \rank H^0(\Gamma,\Z)=1 $, and we have
\begin{align*}
g(\Gamma)=1 - \chi(\Gamma) = 1-|V(\Gamma)|+|E(\Gamma)|,
\end{align*}
where $ \chi(\Gamma) =\rank H^0(\Gamma,\Z) - \rank H^1(\Gamma,\Z)  $ is the Euler characteristic of $ \Gamma $. 

A subgraph $ H $ of $ \Gamma $ is a finite union of edges (possibly disconnected), and $ \Gamma-H $ means the closure of $ \{x\in \Gamma:x\notin H\} $ in $ \Gamma $, hence $ \Gamma - H $ is also a subgraph of $ \Gamma $, and $ H\cap (\Gamma-H) $ is a finite union of vertices. 

Let $ \mathrm{Zh}(\Gamma) $ be the set of all continuous functions $ f:\Gamma\to \R $ such that $ f $ is $ C^2 $ on $ \Gamma\backslash V(\Gamma) $ and $ f''(x)\in L^1(\Gamma) $. For a function $ f\in\mathrm{Zh}(\Gamma)  $, Chinburg and Rumely \cite{CR} defined a measure-valued Laplacian operator by 
\begin{align*}
\Delta f = -f''(x) d x -\sum_{p\in V(\Gamma)} \left(\sum_{\vec{v} \text{ at } p} d_{\vec{v}} f(p)  \right) \delta_p(x).
\end{align*}
Here $ dx $ is the standard Lebesgue measure on $ \Gamma $, $ \delta_p(x) $ is the Dirac measure at $ p $, $ \vec{v} $ is a unit tangent vector at $ p $, and $ d_{\vec{v}} f(p)=\lim_{t\to 0^+} \frac{f(p+t\vec{v})-f(p)}{t} $. 

We list some basic properties that will be used later. For their proofs, see \cite{CR,Zha93}. 

\begin{prop}\label{Delta}
	(1) (Self-adjointness) For each $ f,g\in \mathrm{Zh}(\Gamma) $, we have
	\begin{align*}
	\int_\Gamma g  \Delta f=\int_\Gamma f \Delta g=\int_\Gamma f'(x) g'(x) dx.
	\end{align*}
	As a corollary, for each $ f\in \mathrm{Zh}(\Gamma) $, 
	\begin{align*}
	\int_\Gamma  \Delta f=0.
	\end{align*}
	i.e., $ \Delta f $ has total mass zero. Further, $ \Delta f = 0 $ if and only if $ f $ is a constant. 
	
	(2) (Additivity) If $ \Gamma $ is the union of two subgraphs $ \Gamma_1 $ and $ \Gamma_2 $, such that $ \Gamma_1\cap \Gamma_2 $ is a finite set of points. Then we have
	\begin{align*}
	\Delta f=\Delta (f|_{\Gamma_1}) + \Delta (f|_{\Gamma_2}),
	\end{align*}
	where $ \Delta (f|_{\Gamma_i}) $ is a measure on $ \Gamma_i $, which is also a measure on $ \Gamma $ by pushforward of measures. 
	
	(3) (Maximal principle) If $ \Delta f $ is a finite union of dirac measures, say,
	\begin{align*}
	\Delta f=\sum_{i=1}^n c_i \delta_{x_i}(x),
	\end{align*} 
	where $ c_i\neq 0 $ for all $ i $. Then $ f $ achieves maximum (resp. minimum) at some $ x_i $ with $ c_i > 0 $ (resp. $ c_i< 0 $). 
\end{prop}
Note that $ \Delta(f|_{\Gamma_1}) $ may not equal to $ (\Delta f)|_{\Gamma_1} $, and there exists real numbers $ c_i $ such that
\begin{align*}
\Delta(f|_{\Gamma_1})-(\Delta f)|_{\Gamma_1} = \sum_{x_i\in \Gamma_1\cap (\Gamma- \Gamma_1)} c_i \delta_{x_i}. 
\end{align*}

In \cite[\S 2]{CR}, they introduced the voltage function $ j_z(x,y) $. 

\begin{prop}
	Fix $ y,z\in \Gamma $. Then there exists a unique function $ j_z(x,y) $ on $ \Gamma $, continuous and piecewise linear in $ x $, such that 
	\begin{align*}
	\begin{cases}
	\Delta_x j_z(x,y)=\delta_y(x) - \delta_z (x), \\
	j_z(z,y) = 0. 
	\end{cases}
	\end{align*}
	The function $ j_z(x,y) $ is symmetric in $ x $ and $ y $, and is jointly continuous in all three variables $ x,y,z $. 

\end{prop}

	Define the resistance function $ r(x,y)=j_y(x,x) $ for $ x,y\in \Gamma $ as in \cite{BF,BR,Cin}. If we want to emphasize the graph, we will write $ j_z(\Gamma;x,y) $ and $ r(\Gamma;x,y) $. 
	
	Then by Proposition \ref{Delta}, we get
	\begin{align*}
	0=j_z(z,y)\leqslant j_z(x,y)\leqslant j_z(y,y)=r(y,z) 
	\end{align*}
	for all $ x\in \Gamma $, and $ r(y,z)=j_y(x,z)+j_z(x,y) $ for all $ x\in \Gamma $, thus $ r(y,z) $ is symmetric in $ y,z $.

For any real-valued, signed Borel measure $ \mu $ on $ \Gamma $, such that $ \mu(\Gamma)=1 $ and $ |\mu |(\Gamma)<\infty $, let 
$$ j_\mu(x,y)=\int_\Gamma j_z(x,y) d\mu(z). $$ 
Then the constant
\begin{align*}
 c_\mu(\Gamma) = \int_\Gamma j_\mu(x,y)d\mu(x) 
\end{align*} 
is independent of $ y $, see \cite[Lem. 2.16]{CR}. 

The Green's function $ g_\mu(x,y) $ on $ \Gamma\times \Gamma $ associated to $ \mu $ is defined by \begin{align*}
 g_\mu(x,y)=j_\mu(x,y)-c_\mu(\Gamma). 
\end{align*}
It is symmetric, continuous, and for each $ y\in \Gamma $, we have
\begin{align*}
&\Delta_x g_\mu(x,y) = \delta_y(x)-\mu,\\
&\int_\Gamma g_\mu(x,y) d\mu(x)  = 0.
\end{align*}
Further, $ g_\mu(x,y) $ is characterized by these two properties, see \cite{CR,BR,Zha93} for details.

There are two measures of particular interest. One is the canonical measure $ \mu_{\can} $ defined in \cite[Thm. 2.11]{CR}. It is the unique measure such that $ j_{\mu_{\can}} (x,x) $ is a constant. This measure has the form
\begin{align*}
\mu_{\can} = \sum_{p\in V(\Gamma)} \left( 1-\frac{1}{2} v(p) \right)\delta_p(x) + 
\sum_{e\in E(\Gamma)} \frac{dx}{ L(e) + R(e)},
\end{align*}
where $ dx $ is the standard Lebesgue measure on the edge $ e $, $ L(e) $ is the length of $ e $, and $ R(e)= r(\Gamma-e; p,q)\in \R\cup \{\infty\} $, where $ p,q $ are the endpoints of the edge $ e $. 

The other one is Zhang's admissible measure defined in \cite{Zha93}, which needs a polarization on $ \Gamma $. 

\subsubsection*{Polarization on Graphs}

Let $ \Gamma $ be a graph. A polarization on $ \Gamma $ is a function $ q:V(\Gamma)\to \Z_{\geqslant 0} $, such that for each $ p\in V $, $ v(p)-2+2q(p)\geqslant 0 $. The pair $ (\Gamma, q) $ is called a polarized graph. 

The formal sum
\begin{align*}
K=\sum_{p\in V(\Gamma)} \left(v(p)-2+2q(p) \right) p
\end{align*}
is called the canonical divisor on $ (\Gamma,q) $, and it defines a positive measure 
\begin{align*}
\delta_K = \sum_{p\in V(\Gamma)} (v(p)-2+2q(p)) \delta_p.
\end{align*}
The degree of the divisor $ K $ is
\begin{align*}
\deg K = \sum_{p\in V(\Gamma)} \left(v(p)-2+2q(p) \right) = \delta_K(\Gamma). 
\end{align*}

The genus of a polarized graph $ (\Gamma,q) $ is defined to be
\begin{align*}
g(\Gamma,q) = g(\Gamma)+\sum_{p\in V(\Gamma)} q(p), 
\end{align*}
and we see that $ \deg K = 2g(\Gamma,q)-2 $. 

The admissible measure $ \mu_{\ad} $ is defined in \cite{Zha93}. It is the unique measure such that the function $ g_{\mu_{\ad}} (x,x)+ g_{\mu_{\ad}}(x,\delta_K) $ is a constant. It has the form
\begin{align*}
\mu_{\ad} = \frac{1}{2g} \left( 2\mu_{\can} + \delta_K \right).
\end{align*}

The corresponding function $ g_{\mu_{\mathrm{ad}}}(x,y) $ is called the admissible Green's function. 

Many invariants coming from algebraic geometry can be deduced from $ \delta_K $, $ \mu_{\can} $, $ g_{\mu}(x,y) $, and $ r(x,y) $.

Zhang's $ \varepsilon $-invariant in \cite{Zha93} is defined by
\begin{align*}
\varepsilon(\Gamma,q) &= \int\int _{\Gamma\times \Gamma} r(x,y)\delta_K(x)\mu_{\ad}(y).
\end{align*}
Zhang's $ \phi $-invariant in \cite{Zha10} is defined by
\begin{align*}
\phi(\Gamma,q) &= -\frac{1}{4} \ell(\Gamma) + \frac{1}{4} \int_{\Gamma} g_{\mu_{\ad}}(x,x) ((10g+2)\mu_{\ad} -\delta_{K}).
\end{align*}
Other invariants can be found in \cite{Cin}. We mainly consider the $ \phi $-invariant in this paper.

\subsection{Degeneration of Graphs}\label{DegenerationofGraphs}
In this section, let $ (\Gamma,q) $ be a polarized graph. We first define a notion of degeneration of graphs, and then we show that most invariants and functions discussed in previous section are continuous under degeneration. 

First we omit the polarization. Let $ \Gamma $ be a fixed graph. Then for each functions $ L:E(\Gamma)\to \R_{>0} $, we get a metrized graph $ (\Gamma,L) $, such that the length of the edge $ e\in E(\Gamma) $ is $ L(e) $. All metrized graphs with underlying graph $ \Gamma $ can be obtained this way. We write $ L\in \R_{>0}^{E(\Gamma)} $ to denote the function. In other word, $ \R_{>0}^{E(\Gamma)} $ is the "moduli space" of metrized graphs with the given underlying graph $ \Gamma $.

Now $ \R_{>0}^{E(\Gamma)} $ is a subset of $ \R^{E(\Gamma)} $, where the latter space is endowed with a standard Euclidean metric. Then  $ \R_{\geqslant0}^{E(\Gamma)} $ is the closure of $ \R_{> 0}^{E(\Gamma)} $. We will give a moduli interpretation of $ \R_{\geqslant 0}^{E(\Gamma)} $. 

\begin{definition}
	Let $ L\in \R_{\geqslant 0}^{E(\Gamma)} $. The metrized graph $ (\Gamma,L) $ is defined as follows. Let 
	$$ E_0=\{ e\in E(\Gamma): L(e)=0 \}, $$ 
	viewed as a subgraph of $ \Gamma $. The underlying topological space of $ (\Gamma,L) $ is the quotient space $ \Gamma/E_0 $. Then $ E(\Gamma/E_0)=E(\Gamma)\backslash E_0 $, and the metric on $ \Gamma/E_0 $ is given by the restriction $ L|_{E(\Gamma/E_0)}:E(\Gamma/E_0)\to \R_{>0} $. 
\end{definition} 

Let $ L\in \R_{\geqslant 0}^{E(\Gamma)} $. The vertex set of $ (\Gamma,L) $ is a quotient of $ V(\Gamma) $, and the edge set of $ (\Gamma,L) $ is a subset of $ E(\Gamma) $. 

Now consider the polarization. Let $ (\Gamma,q) $ be a polarized graph.

\begin{definition}\label{Def2.4}
	Let $ L\in \R_{\geqslant 0} ^{E(\Gamma)} $. Let $ F:\Gamma\to (\Gamma,L) $ be the canonical quotient map. Let $ p' $ be any vertex of $ (\Gamma,L) $. Define the polarization $ q' $ on $ (\Gamma,L) $ as
	\begin{align*}
	q'(p')=\sum_{p\in F^{-1}(p')} q(p)+g(F^{-1}(p'))
	\end{align*}
	where $ F^{-1}(p') $ is a connected subgraph of $ \Gamma $, and $ g(F^{-1}(p')) $ is its genus.
\end{definition}

We need to check that $ q' $ is indeed a polarization on $ (\Gamma,L) $. 
\begin{prop}
	The metrized graph $ (\Gamma,L) $ together with the function $ q' $ is a polarized metrized graph with genus $ g((\Gamma,L),q')=g(\Gamma,q) $. 
\end{prop}

\begin{proof}
	We first show that $ \delta_{K'}=F_* \delta_{K} $, where $ K $ is the canonical divisor on $ (\Gamma,q) $, and $ K' $ is the canonical divisor of $ ((\Gamma,L),q') $. Write $ F_*\delta_K=\sum a_i p_i $, and let $ E_i = F^{-1}(p_i) $. Then we get
	\begin{align*}
	a_i&=\delta_K(E_i) = \sum_{p\in V(E_i)}\left( v(p) -2 + 2q(p)\right) \\
	&  =v(p_i)+2|E(E_i)|-2|V(E_i)| + 2\sum_{p\in V(E_i)}q(p)\\
	&=v(p_i) -2 + 2g(E_i)+ 2\sum_{p\in V(E_i)}q(p)\\
	&= v(p_i)-2+2q'(p_i)=\delta_{K'}(p_i). 
	\end{align*}
	Hence $ \delta_{K'}=F_* \delta_K $ is a positive measure, and $ ((\Gamma,L),q') $ is a polarized metrized graph with genus $ g((\Gamma,L),q')=1+\frac{1}{2} \deg K' =1+\frac{1}{2}\deg K =g(\Gamma,q) $. 
\end{proof}
By abuse of notation, we use $ (\Gamma,L,q) $ to denote this polarized graph, although the polarization is actually $ q' $ instead of $ q $. Note that even if $ q=0 $, the polarization $ q' $ on $ (\Gamma,L,q) $ may be nonzero.

In both cases, we get a moduli interpretation of $ \R_{\geqslant0}^{E(\Gamma)} $.

\subsection{Continuity Properties}\label{ContinuityProperties}
Now we consider the continuity properties of (polarized) metrized graph $ (\Gamma,L) $ (or $ (\Gamma,L,q)$) with respect to $ L $. We ignore the polarization to simplify the notations.

We choose a particular metric 
$$ L_0:E(\Gamma)\to \R,\ \  L_0(e)=1 \text{ for all } e\in E(\Gamma), $$
and we get a particular metrized graph $ (\Gamma,L_0) $. Then for each $ L\in \R_{\geqslant 0}^{E(\Gamma)} $, there is a unique piecewise linear morphism $ F_{L}: (\Gamma,L_0)\to (\Gamma,L) $ which is either the identity or the quotient map of graphs.

Let $ n\geqslant 0 $, and $ f:(\Gamma,L)^n\to \R $ be a continuous function on $ (\Gamma,L)^n $. Its pull-back $ F_L^*f:(\Gamma,L_0)^n\to \R $ is given by 
\begin{align*}
F_L^* f(x_1,\ldots,x_n)= f(F_L(x_1),\ldots,F_L(x_n)),
\end{align*}
and if $ n=0 $, the function $ f $ is just a real number, and we don't really need the pull-back map $ F_L^* $. 

Now for different $ L $, functions on $ (\Gamma,L)^n $ can be all viewed as functions on $ (\Gamma,L_0) $, so they can be compared with others. 

Suppose that for each $ L\in \R^{E(\Gamma)}_{>0} $, a function $ f_L:(\Gamma,L)^n \to \R $ is given, which is denoted by $ f(\Gamma;L;x_1,\ldots,x_n) $. Be careful that for different $ L $, the variables $ x_i $ lie in different space.

\begin{definition}
	The function $ f(\Gamma;L;x_1,\ldots,x_n) $ is called continuous under degeneration of graphs if the following holds.
	
	For any sequence of points $ L_i\in \R_{\geqslant 0}^{E(\Gamma)} $ converging to $ \tilde{L}\in \R_{\geqslant 0}^{E(\Gamma)} $, the sequence of functions $ F_{L_i}^* f(\Gamma; L_i;x_1,\ldots,x_n) $ converges to $ F^*_{\tilde{L}} f(\Gamma;\tilde{L};x_1,\ldots,x_n) $ uniformly on $ \Gamma^n $. 
\end{definition}
Since $ \R_{>0}^{E(\Gamma)} $ is dense in $ \R_{\geqslant 0}^{E(\Gamma)} $, to check the continuity, it suffices to check all sequences $ L_i\in \R^{E(\Gamma)}_{>0} $ converging to $ \tilde{L}\in \R_{\geqslant 0}^{E(\Gamma)} $. 

We will show that most invariants and functions are continuous under degeneration of graphs. 

\subsubsection*{Continuity of the Voltage Function}

We first show that the voltage function $ j(\Gamma;L;x,y,z)=j_z(\Gamma;L;x,y) $ is continuous under degeneration of graphs. Before the proof, we do some preparations. Let $ (\Gamma,L) $ be a fixed metrized graph.

\begin{lemma}\label{Slope}
	 For each fixed $ y,z\in \Gamma $, and for each tangent direction $ v $ at $ x $, we have $ |d_{v,x}(j(x,y,z))|\leqslant 1 $. 
\end{lemma}

\begin{proof}
	By adding points, we may assume that $ x,y,z $ are all vertices in $ \Gamma $. 
	
	Let $ t=j(x,y,z) $, and let $ \Gamma_t=\{x\in \Gamma: j_z(x,y)\leqslant t\} $. 
	
	If $ x\neq y $, then $ \Gamma_t\neq \Gamma $ and let $ \Gamma'_t=\Gamma-\Gamma_t $. Then $ \Gamma_t\cap \Gamma'_t $ is a finite set, say, $ x_1,\ldots,x_n $. If $ x=y $, let $ n=1 $ and $ x_1=y $.

	By additivity of laplacian, we have that 
	\begin{align}\label{Lap}
	\Delta_p(j(p,y,z)|_{\Gamma_t}) = -\delta_z(p)+\sum_{i=1}^n c_i \delta_{x_i}(p).
	\end{align}
	for some real numbers $ c_i $. Since the Laplacian has total volume $ 0 $, $ \sum_{i=1}^n c_i= 1 $. 
	
	Now, since $ j|_{\Gamma_t} $ achieves maximum at $ x_i $, for each inner direction $ v $ of $ \Gamma_t $ at $ x_i $, we have $ d_{v}(j(x_i,y,z))\leqslant 0 $, and hence $ c_i=-\sum_{v} d_{v} j(x_i,y,z)\geqslant 0 $ and thus $ c_i\leqslant 1 $ and thus $ -1 \leqslant d_{v} j(x_i,y,z)\leqslant 0 $. 
	
	If $ x = x_i $ for some $ i $, and $ v $ is a inner direction at $ x $ of $ \Gamma_t $, we've done. If $ v $ is an inner direction of $ \Gamma_t' $, similarly we get $ 0\leqslant d_{v,x} j(x,y,z)\leqslant 1 $. In both cases, we have $ |d_{v,x} j(x,y,z)|\leqslant 1 $. 
	
	Now assume $ x\notin \Gamma'_t $. Since $ j|_{\Gamma_t} $ also achieves maximum at $ x $, it is constant by the equation  (\ref{Lap}). Thus $ d_{v,x}(j(x,y,z))=0 $ for all tangent direction $ v $ at $ x $. 
	\end{proof}

Let $ (\Gamma,L) $ be a polarized graph with $ L\in \R_{>0}^{E(\Gamma)} $. 	For each edge $ e $ in $ \Gamma $, we give an orientation on $ e $. More precisely, let $ p,q $ be its two endpoints(possibly equal). Choose one endpoint, say $ p $, as the source, and the other point $ q $ is the target of $ e $. We write $ p=s(e) $ and $ q=t(e) $. 

We choose two vertices $ y\neq z $ in $ \Gamma $. Let $ \{x_e\}_{e\in E(\Gamma)} $ be indeterminants. Consider the following system of linear equations.

For each vertex $ v\neq z $ in $ \Gamma $, there is an equation
\begin{align}\label{Flow1}
	\sum_{e:v=t(e)}  x_{e} - \sum_{e:v=s(e)} x_{e} = \varepsilon_{vy},
\end{align}
where $ \varepsilon_{vy}=1 $ if $ v=y $, and otherwise $ \varepsilon_{vy}=0 $. 

Further, let $ \{\gamma_1,\ldots,\gamma_g\} $ be loops in $ \Gamma $ such that they become a basis of the homology group $ H_1(\Gamma,\Z) $. For each $ \gamma_j $, there is an equation
\begin{align}\label{Flow2}
\sum_{e\in \gamma_j} \pm L(e) x_{e} = 0,
\end{align}
where the sign depends on whether the orientation of $ e $ agrees with the orientation of $ \gamma_j $. 

\begin{lemma}\label{LemmaofEquations}
	The system of equations (\ref{Flow1}) and (\ref{Flow2}) has a unique soluation \begin{align}\label{Solution}
	x_e=\frac{1}{L(e)}\left(j(t(e),y,z)-j(s(e),y,z)\right).
	\end{align}
\end{lemma}
\begin{proof}
	We note that there are in total $ |V(\Gamma)|-1+g(\Gamma)=|E(\Gamma)| $ equations and indeterminants.
	
	First we check that the formula (\ref{Solution}) is indeed a solution of the equations. Just note that the equation (\ref{Flow1}) is a rewrite of the Laplacian equation, and the equation (\ref{Flow2}) is clear. 
	
	Now we check the uniqueness. Write equations (\ref{Flow1}) and (\ref{Flow2}) as $ Ax = b_y $. Suppose that $ A\alpha=0 $ for some $ \alpha=(a_e) $. Fix a vertex $ v\in V(\Gamma) $. Then for each other vertices $ w\in V(\Gamma) $, choose a path $ \gamma=[v=v_0,e_1,v_1,e_2,\ldots e_n,v_n=w] $ connecting $ v $ ans $ w $. Define 
	\begin{align*}
	f(w)= \sum_{i=1}^n \varepsilon_i L(e_i) a_{e_i}
	\end{align*}
	where $ \varepsilon_i=1 $ if $ v_i $ is the target endpoint of $ e_i $, and otherwise $ \varepsilon_i=-1 $. Note that $ f(w) $ is independent of the path due to equation (\ref{Flow2}).
	
	Now we extends $ f $ piecewise linearly. Then $ \Delta f=0 $, with $ f(v)=0 $. Hence $ f=0 $, and $ a_{e_i}=(f(t(e_i))-f(s(e_i)))/L(e_i)=0 $. Thus the solution is unique.
\end{proof}

\begin{theorem}\label{Voltage}
	The voltage function $ j(\Gamma;L;x,y,z)=j_z(\Gamma;L;x,y) $ is continuous under degeneration of graphs. 
\end{theorem}

\begin{proof}
	Take any sequence $ L_m\in \R_{>0}^{E(\Gamma)} $ converging to $ \tilde{L}\in \R_{\geqslant 0}^{E(\Gamma)} $. We ignore the pull-back notation $ F_{L_m}^* $  to simplify the notations, but we still write $ \tilde{F} $ to emphasis that $ (\Gamma,\tilde{L}) $ may not be homeomorphic to $ (\Gamma,L_0) $. 
	
	Let $ j_m(x,y,z) $ to be the pull back of $ j_{L_m}(x,y,z) $ on $ (\Gamma,L_0) $, and let $ \Delta_m $ be the pull-back of Laplacian operator of $ (\Gamma,L_i) $ on $ (\Gamma,L_0) $. 
	
	First, we fixed $ y,z $, and consider the continuity property of $ j(x,y,z) $ in $ x $. Assume that $ y,z $ are both vertices on $ \Gamma $. If $ y=z $, then $ j(\Gamma;L;x,y,z)=0 $ for all $ L $. Hence we assume $ y\neq z $
	
	Fix an orientation for each vertex $ v\in V(\Gamma) $, and choose loops $ \gamma_1,\ldots,\gamma_g $ of $ \Gamma $ such that they become a basis of $ H^1(\Gamma,\Z) $. We may assume that the first $ s $ loops $ \gamma_1,\ldots,\gamma_s $ maps to a basis of $ H_1((\Gamma,\tilde{L}),\mathbb{Z}) $.
	
	Let $ a_{e,m}=(j_m(t(e),y,z)-j_m(s(e),y,z))/L_m(e) $. Then $ a_{e,m} $ satisfies the following system of linear equations for the metrized graph $ (\Gamma,L_m) $ by Lemma \ref{LemmaofEquations}. 
	\begin{align*}
	\begin{cases}
		\ds\sum_{e:v=t(e)}  x_{e} - \ds\sum_{e:v=s(e)} x_{e} = \varepsilon_{vy},\quad&\text{where }\  v\in V(\Gamma)\ \text{ and }v\neq z,\\
		\ds\sum_{e\in \gamma_j} \pm L_m(e) x_{e} = 0, \qquad \qquad& 1\leqslant j\leqslant g,
	\end{cases}
	\end{align*}
	where $ \varepsilon_{vy}=1 $ if $ v=y $, and otherwise $ \varepsilon_{vy}=0 $. 
	
	Next, we change the equations so that the limit of the equations is exactly the corresponding equation of $ (\Gamma,\tilde{L}) $. 
	
	Let $ \tilde{F}:\Gamma\to (\Gamma,\tilde{L}) $ be the canonical quotient map. For each point $ p\in V(\Gamma,\tilde{L}) $ with $ p\neq \tilde{F}(z) $, let $ E_p=\tilde{F}^{-1}(p) $ be the preimage subgraph. Then $ \{a_{e,m}\} $ satisfies the equation
	\begin{align}\label{Eq1}
	\sum_{v\in E_p} \left(\sum_{e:v=t(e)}  a_{e,m}-\sum_{e:v=s(e)} a_{e,m} \right) = \varepsilon_{py},
	\end{align}
	where $ \varepsilon_{py}=1 $ if $ y\in E_p $, and otherwise $ \varepsilon_{py}=0 $. Then if $ e\in E_p $, the term $ a_{e,m} $ does not occur in the above equation.

	For each loop $ \gamma_j $ with $ 1\leqslant j\leqslant s $, we have
	\begin{align*}
	\sum_{e\in\gamma_j,e\notin E(\tilde{L})} \pm L_m(e) a_{e,m} = 		-\sum_{e\in\gamma_j,e\in E(\tilde{L})} \pm L_m(e) a_{e,m}.
	\end{align*}
	
	For each $ m $, consider the following system of linear equations of indeterminants $ \{x_e\}_{e\notin E(\tilde{L})} $
\begin{equation}
	\begin{cases}
	\displaystyle\sum_{v\in E_p} \left(\displaystyle\sum_{e:v=s(e)}  x_e-\displaystyle\sum_{e:v=t(e)} x_e \right) & = \varepsilon_{py}, \\
	\displaystyle\sum_{e\in\gamma_j,e\notin E(M)} \pm L_m(e) x_e & = 		-\displaystyle\sum_{e\in\gamma_j,e\in E(M)} \pm L_m(e) a_{e,m}
	\end{cases}
\end{equation}
where $ p\neq \tilde{F}(z) $ is a vertex in $ (\Gamma,\tilde{L}) $, and $ 1\leqslant j\leqslant s $. We denote it by $ A_m x = b_m $. Then $ \alpha_m= \{a_{m,e}\}_{e\in E(\Gamma)\backslash E(M)} $ is a solution of $ A_m x=b_m $. 
	
	Now let $ m\to \infty $. The limit $ A= \displaystyle\lim_{m\to \infty} A_m $ and $ b=\displaystyle\lim_{m\to\infty } b_m $ exists, and the equation $ A x=b $ is exactly the equations of $ (\Gamma,\tilde{L}) $ as in Lemma \ref{LemmaofEquations}. Thus $ \det A\neq 0 $ and $ A $ is invertible, which implies
	\begin{align*}
	\lim_{m\to \infty} A_m^{-1} = A^{-1}
	\end{align*} 
	by the continuity of the inverse of matrix.
	
	Now let
	\begin{align*}
	 \tilde{\alpha} = \lim_{m\to\infty} \alpha_m=\lim_{m\to\infty}A_m^{-1} b_m=  A^{-1} b,
	\end{align*}
	then $ \tilde{\alpha} $ satisfies $ A\tilde{\alpha}=b $.

	Now let's go back to $ j_z(x,y) $. The above proposition shows that for each $ e\notin E(\tilde{L}) $, we have 
	$$ \lim_{m\to \infty} \frac{j_z(s(e),y)-j_z(t(e),y)}{L_m(e)} =\frac{ \tilde{F}^* j_{z}(s(e),y)-\tilde{F}^* j_{z}(t(e),y)}{\tilde{L}(e)}. $$
	As $ \tilde{L}(e)\neq 0 $, we have
	$$ \lim_{m\to \infty} j_z(s(e),y)-j_z(t(e),y) =  \tilde{F}^* j_{z}(s(e),y)-\tilde{F}^* j_{z}(t(e),y)). $$
	
	Further, if $ e\in E_0 $, then $ |j_m(s(e),y,z)-j_m(t(e),y,z)|\leqslant L_{m}(e) $ by Lemma \ref{Slope}, hence we get
	\begin{align*}
	\lim_{m\to \infty} j_z(s(e),y)-j_z(t(e),y) = 0 = \tilde{j}_z(s(e),y)-\tilde{j}_z(t(e),y).
	\end{align*}
	As $ j_m(z,y,z)=\tilde{F}^* j(z,y,z)=0 $ for all $ m $, we see $ j_m(x,y,z) $ converges to $ \tilde{F}^* j(x,y,z) $ for all vertices $ x $. As $ j_m(x,y,z),\tilde{F}^* j(x,y,z) $ are all piecewise linear in $ x $, we show that $ j_m(x,y,z) $ converges to $ \tilde{F}^* j(x,y,z) $ uniformly in $ x $. 
	
	Now we consider the uniformity in $ y,z $. Since $ \Gamma $ is compact, we only need to show the uniformity on a closed interval. Fix $ y_0,z_0 $. Given any $ \varepsilon>0 $, then choose $ N>0 $ such that $ |j_n(x,y_0,z_0)-\tilde{j}(x,y_0,z_0)|<\varepsilon $ for all $ x $. Then for all $ y,z $ with $ |y-y_0|\leqslant \varepsilon $, $ |z-z_0|\leqslant \varepsilon $, we have	
	\begin{align*}
	&|j_n(x,y,z)-\tilde{j}(x,y,z)|\\ & \leqslant |j_n(x,y,z)-j_n(x,y_0,z_0)|+|j_n(x,y_0,z_0)-\tilde{j}(x,y_0,z_0)| \\
	 & \qquad + |\tilde{j}(x,y_0,z_0)-\tilde{j}(x,y,z)|\\
	&\leqslant 2|y-y_0|+2|z-z_0|+|j_n(x,y_0,z_0)-\tilde{j}(x,y_0,z_0)|\leqslant 5\varepsilon. 
	\end{align*}
	Thus the converges is uniformly in all $ x,y,z $. 
\end{proof}

\subsubsection*{Continuity  of Measures}
Next, we discuss the continuity property of measures under the degeneration of graphs. 

\begin{definition}
	Let $ \Gamma $ be a graph with or without polarization. Suppose that for each $ L\in \R_{\geqslant 0}^{E(\Gamma)} $, there is a measure $ \mu(\Gamma,L) $ on $ (\Gamma,L) $. We say that the measures are continuous under degeneration if the followings hold:
	
	(1) For each $ L_m\in \R_{>0}^{E(\Gamma)} $ converging to $ \tilde{L} $, the norm of the measures $ |\mu(\Gamma,L_m)| $ is bounded. 
	
	(2) We view $ \mu(\Gamma,L_m) $ as measures on $ (\Gamma,L_0) $ via pushforward of measures by $ F_m^{-1} $. Let $ \tilde{F}:(\Gamma,L_0)\to (\Gamma,\tilde{L}) $ as before. Then we have
	\begin{align*}
	|\tilde{F}_*\mu(\Gamma,L_m) - \mu(\Gamma,\tilde{L})|\to 0
	\end{align*}
\end{definition}
We sill show that the canonical measure and the admissible measure are both continuous under degeneration. 

\subsubsection*{Continuity of the Canonical Measure}

Let $ (\Gamma,L) $ be a metrized graph. The canonical measure is defined as
\begin{align*}
\mu_{\can}(x) = \sum_{p\in V(\Gamma)} \left(1-\frac{1}{2}v(p) \right)\delta_p(x)+\sum_{e\in E(\Gamma)} \frac{dx}{L(e)+R(e)}.
\end{align*}
Here $ R(e) $ is defined as follows. Let $ e^0 $ be the interior of the edge, and let $ p,q $ be the endpoints of $ e $. If $ \Gamma-e^0 $ is disconnected, then $ R(e)=\infty $; otherwise, 
\begin{align*}
R(e)=r(\Gamma-e^0,L;p,q).
\end{align*}

Chiburg and Rumely \cite[\S 2]{CR} showed that for any metrized graph $ \Gamma $, 
\begin{align*}
\sum_{e\in E(\Gamma)} \frac{L(e)}{L(e)+R(e)} = g(\Gamma).
\end{align*}

To prove the continuity of the canonical measure, we meed a lemma.

\begin{lemma}\label{CurrentLemma}
	Let $ (\Gamma,L) $ be a metrized graph. Fix $ 0 < \varepsilon < A $. Let $ E_0 $ be a connected subgraph of $ \Gamma $. Suppose that $ L(e)\leqslant\varepsilon $ for all $ e\in E_0 $, and $ L(e)\geqslant A $ for any edge $ e\notin E_0 $ intersecting $ E_0 $. Then for any two vertices $ p,q\in E_0 $, we have
	\begin{align*}
	\left|\frac{r(\Gamma;p,q)-r(E_0;p,q)}{r(\Gamma;p,q)}\right| \leqslant {2|E(\Gamma)|^2} \frac{\varepsilon}{A}
	\end{align*}
\end{lemma}

\begin{proof}
	Let's consider the resitriction of $ j(\Gamma;x,p,q) $ on $ E_0 $. By the additivity of Laplacian, we have
	\begin{align*}
	\Delta_x \left(j(\Gamma;x,p,q)|_{E_0}\right) = \delta_p-\delta_q+\sum_{i} c_{i} \delta_{x_i}(x),
	\end{align*}
	where $ x_i\in E_0\cap (\Gamma\backslash E_0)  $, and $ c_i $ are real numbers.
	
	Note that $ r(\Gamma;x,y)\leqslant | E(\Gamma)| \varepsilon $ for all $ x,y\in E_0 $. In particular, $ r(\Gamma;p,q)\leqslant | E(\Gamma)| \varepsilon $. Thus for any two vertices $ x',x''\in \Gamma $ on the same edge, we have
	\begin{align*}
	|j(\Gamma;x',p,q)-j(\Gamma;x'',p,q)|\leqslant r(\Gamma;p,q)\leqslant |E(\Gamma)| \varepsilon.
	\end{align*}
	
	Note that $ c_i $ can also be computed in $ \Gamma-E_0 $, and since the length of any edges intersecting $ E_0 $ has a uniform positive lower bound $ A $, we get
	$$ |c_i|\leqslant \frac{v(x_i)r(\Gamma;p,q)}{A}. $$
	
	Consider two functions $ j(E_0;x,p,q) $ and $ j(\Gamma;x,p,q)|_{E_0} $ on $ E_0 $. The Laplacian of their difference is
	$$ \Delta_{x} \left( j(\Gamma;x,p,q)|_{E_0} -  j(E_0;x,p,q) \right)=\sum c_i \delta_{x_i} (x), $$ 
	which is a finite sum of Dirac measures, hence we get
	\begin{align*}
	j(\Gamma;x,p,q)|_{E_0} -  j(E_0;x,p,q)=\sum c_i j(E_0;x,x_i,q),
	\end{align*}
	as they have the same Laplacian and vanish at $ q $. Thus 
	\begin{align*}
	|j(\Gamma;x,p,q)-j(E_0;x,p,q)|&\leqslant \sum_i |c_i| r(E_0;x_i,q)\\
	& \leqslant \sum_i \frac{v(x_i)r(\Gamma;p,q)}{A}|E(\Gamma)| \varepsilon \\& \leqslant \frac{2|E(\Gamma)| r(\Gamma;p,q)}{A} |E(\Gamma)| \varepsilon.
	\end{align*}
	
	Take $ x=p $, we get
	\begin{align*}
	\left|\frac{r(\Gamma;p,q)-r(E_0;p,q)}{r(\Gamma;p,q)}\right| \leqslant {2|E(\Gamma)|^2} \frac{\varepsilon}{A}.
	\end{align*}
\end{proof}

Now we consider the continuity of canonical measures under degeneration of graphs. As before, let $ L_m\in \R_{>0}^{E(\Gamma)} $ converging to $ \tilde{L}\in \R_{\geqslant 0}^{E(\Gamma)} $, and let $ \tilde{F}:(\Gamma,L_0)\to (\Gamma,\tilde{L}) $ be the quotient map. We view $ \mu_{m,\can} $ as a measure on $ (\Gamma,L_0) $ to simplify the notations, and let $ \tilde{\mu}_{\can} $ be the canonical measure on $ (\Gamma,\tilde{L}) $. Clearly that $ |\mu_{m,\can}| $ is uniformly bounded, so it suffices to check the second condition. 

\begin{prop}
	The measures $ \tilde{F}_* \mu_{m,can} $ converges to $ \tilde{\mu}_{\can} $ in the sense that $ |\tilde{F}_* \mu_{m,can}-\tilde{\mu}_{\can}|\to 0 \ $. 
\end{prop}

\begin{proof}
	Let $ E_0=\{ e\in E(\Gamma): \tilde{L}(e)=0 \} $. By direct computation, we get
	\begin{align*}
	&\quad \tilde{F}_* \mu_{m,\can}-\tilde{\mu}_{\can}\\
	&=\sum_{p'\in V(\Gamma,\tilde{L})}\left( \mu_{m,\can}(\tilde{F}^{-1}(p'))-\tilde{\mu}_{\can}(p') \right)\delta_{p'}(x) + \sum_{e\in E(\Gamma,\tilde{L})}\left( \mu_{m,\can}(e^0)-\tilde{\mu}_{\can}(e^0)\right)\frac{dx}{\tilde{L}(e)} \\	
	&= \sum_{p'\in V(\Gamma,\tilde{L})}\left(\sum_{\substack{p\in V(\Gamma)\\\tilde{F}(p)=p'}}\left( 1-\frac{v(p)}{2}\right)+\sum_{e\in\tilde{F}^{-1}(p')} \frac{L_m(e)}{L_m(e)+R_m(e)}-1+\frac{v(p')}{2} \right)\delta_{p'}(x)\\
	&\qquad +\sum_{e\in E(\Gamma,\tilde{L})} \left(\frac{L_m(e)}{L_m(e)+R_m(e)}- \frac{\tilde{L}(e)}{\tilde{L}(e)+\tilde{R}(e)}\right)\frac{dx}{\tilde{L}(e)}\\
	&=: I_{1,m}+ I_{2,m}	
	\end{align*}
	where $ dx $ is the Lebesgue measure on $ (\Gamma,\tilde{L}) $, and $ e^0 $ is the interior of $ e $.

	First we consider $ I_{2,m} $. Let $ e\in E(\Gamma,\tilde{L}) $, then $ e $ is also an edge in $ (\Gamma,L_m) $. If $ \Gamma\backslash e^0 $ is disconnected, then so is $ (\Gamma,\tilde{L})\backslash e^0 $, and hence 
	\begin{align*}
	\mu_{m,\can}(e^0)=0=\tilde{\mu}_{\can}(e^0)
	\end{align*}
	Otherwise, we have
	\begin{align*}
	\lim_{m\to \infty} L_m(e)=\tilde{L}(e)> 0,\qquad \lim_{m\to \infty} R_m(e)=\tilde{R}(e)\geqslant 0
	\end{align*}
	by the continuity of the resistance function. Hence $ |I_{2,m}|\to 0 $ as $ m\to\infty $. 
	
	Now we consider $ I_{1,m} $. Let $ p' $ be a vertex in $ (\Gamma,\tilde{L})  $. Let $ E_0=\tilde{F}^{-1}(p') $ be the subgraph. In this case, we compute $ \mu_{m,can} \tilde{F}^{-1}(p) $ as follows.	
	\begin{align*}
	&\ \mu_{m,can} \tilde{F}^{-1}(p') = \sum_{p\in V(E_0)}\left(1-\frac{1}{2}v(p) \right) + \sum_{e\in E_0} \frac{L_m(e)}{L_m(e)+R(\Gamma,L_m,e)}\\
	&=|V(E_0)|- |E(E_0)|-\frac{1}{2}v(p') +\sum_{e\in E_0}\frac{L_m(e)}{L_m(e)+R(\Gamma,L_m,e)}\\
	&=1-\frac{1}{2}v(p')-g(E_0) +\sum_{e\in E_0}  \frac{L_m(e)}{L_m(e)+R(\Gamma,L_m,e)} \\
	\end{align*}
	By the result of Rumely and Baker, we see that
	\begin{align*}
	\sum_{e\in E_0}  \frac{L_m(e)}{L_m(e)+R(E_0,L_m,e)} = g(E_0).
	\end{align*}
	Thus it suffices to prove that 
	\begin{align*}
	\lim_{m\to \infty} \left(\frac{L_m(e)}{L_m(e)+R(E_0,L_m,e)} - \frac{L_m(e)}{L_m(e)+R(E_0,L_m,e)} \right) =0
	\end{align*}
	for all $ e\in E_0 $. 
	
	If $ \Gamma\backslash e^0 $ is disconnected, then $ E_0\backslash e^0 $ is also disconnected, and we have
	\begin{align*}
	\frac{L_m(e)}{L_m(e)+R(\Gamma,L_m,e)} = 0 = \frac{L_m(e)}{L_m(e)+R(E_0,L_m,e)}
	\end{align*}
	
	If $ \Gamma\backslash e^0 $ is connected but $ E_0\backslash e^0 $ is disconnected, then the two end points $ p,q $ of $ e $ map to different points in $ (\Gamma,\tilde{L}) $, thus $ R(\Gamma,L_m,e) $ converges to $ R(\Gamma,\tilde{L},e)> 0 $. Therefore,
	\begin{align*}
	&\quad \lim_{m\to \infty}\quad \left|	\frac{L_m(e)}{L_m(e)+R(\Gamma,L_m,e)} - \frac{L_m(e)}{L_m(e)+R(E_0,L_m,e)} \right|\\
	& =\lim_{m\to \infty}\left|\frac{L_m(e)}{L_m(e)+R(\Gamma,L_m,e)}\right| 
	\end{align*}
	
	If $ \Gamma\backslash e^0  $ and $ E_0\backslash e^0  $ are both connected, then we have
	\begin{align*}
	&\frac{L_m(e)}{L_m(e)+R(\Gamma,L_m,e)} - \frac{L_m(e)}{L_m(e)+R(E_0,L_m,e)}\\
	&=\frac{L_m(e)(R(E_0,L_m,e)- R(\Gamma,L_m,e))}{(L_m(e)+R(\Gamma,L_m,e)) \left(L_m(e)+R(E_0,L_m,e)\right)}.
	\end{align*}
	Note that if $ e'\notin E_0 $ and $ e' $ intersecting $ E_0 $, then  $ L_m(e') $ has a uniform positive lower bound for all $ m $, say, $ L_m(e')\geqslant A>0 $. By the lemma \ref{CurrentLemma} before, we have	
	\begin{align*}
	&\left|\frac{L_m(e)}{L_m(e)+R(\Gamma,L_m,e)} - \frac{L_m(e)}{L_m(e)+R(E_0,L_m,e)}\right| \\	
	&\leqslant \frac{R(\Gamma,L_m,e)}{L_m(e)+R(\Gamma,L_m,e)} \frac{L_m(e)}{L_m(e)+R(E_0,L_m,e)}  \frac{2|E(\Gamma)|^2}{A}\max_{e\in E_0} L_m(e)
	\end{align*}
	Thus $ |I_{1,m}|\to 0 $ as $ m\to \infty $, and $ |\tilde{F}_* \mu_{m,\can}- \tilde{\mu}_{\can}|\to 0 $. 
\end{proof}

\subsubsection*{Other Continuity Properties}
Now we consider the admissible measures. 

Let $ (\Gamma,L,q) $ be a polarized metrized graph with $ L\in \R_{>0}^{E(\Gamma)} $. Let $ \mu_{\ad}(x) $ be the admissible measure associated to $ K $ in \cite{Zha93}. Then
\begin{align*}
\mu_{\ad}(x)=\frac{1}{2g}(2\mu_{\can}(x) + \delta_K(x)).
\end{align*}
Thus $ \mu_{ad} $ is also continuous under degeneration. 

Many invariants of graphs are induced by the voltage function, resistance function, the canonical measure and the admissible measure. Such invariants are continuous under degeneration by the following proposition. 
\begin{prop}
	If the function $ f(\Gamma,L;x_1,\ldots,x_n) $ and the measures $ \mu_{L} $ are continuous under degeneration. Then the function
	\begin{align*}
	g(\Gamma;L;x_2,\ldots,x_n) = \int_{(\Gamma,L)} f(\Gamma,L;x_1\ldots,x_n) d\mu_L(x_1)
	\end{align*}
	is also continuous under degeneration.
\end{prop}
\begin{proof}
	Suppose that $ L_m\to \tilde{L} $ with $ L_m\in \R_{>0}^{E(\Gamma)} $. Let $ \tilde{F}:(\Gamma,L_0)\to (\Gamma,\tilde{L}) $.  Then for each $ x_2,\ldots,x_n\in (\Gamma,L_0) $, we have
	\begin{align*}
	&\left|\int_{(\Gamma,L_m)} f(\Gamma,L_m;x_1\ldots,x_n) d\mu_m(x_1)-\int_{(\Gamma,\tilde{L})} \tilde{F}^* f(\Gamma,\tilde{L};x_1\ldots,x_n) d\tilde{\mu}(\tilde{F}x_1)\right|\\
	& \leqslant \int_{(\Gamma,L_0)} |f(\Gamma,L_m;x_1,\ldots,x_n)- f(\Gamma,\tilde{L};\tilde{F}x_1\ldots,\tilde{F}x_n)| d|\mu_m|(x_1)\\
	&\qquad + \int_{(\Gamma,\tilde{L})} |f(\Gamma,\tilde{L};x_1\ldots,x_n)| d|\tilde{F}_* \mu_{m}-\tilde{\mu}|(x_1)\\
	&\to 0	
	\end{align*}
	Thus the function 
	\begin{align*}
	g(\Gamma;L;x_2,\ldots,x_n) = \int_{(\Gamma,L)} f(\Gamma^{n},L;x_1\ldots,x_n) d\mu(x_1)
	\end{align*}
	is also continuous under degeneration. 
\end{proof}

We know that many functions and invariants are induced by the voltage function, the canonical measure, and the admissible measure. By the above theorem, all these invariants are continuous under degeneration of graphs. 
\begin{corollary}
	The canonical Green function $ g_{\mu_{\can}}(x,y) $, the admissible Green function $ g_{\mu_{\ad}}(x,y) $, the $ \phi $-invariant and the $ \varepsilon $-invariant are all continuous under the degeneration of graphs, 
\end{corollary}
In particular, we get Theorem \ref{PhiforGraph}. 

\section{Admissible Line Bundles}\label{AdmissibleLineBundles}
In this section we will review Yuan--Zhang's work \cite{YZ} on adelic divisors and adelic line bundles, and then review Yuan's work \cite{Yua} on the globalization of admissible line bundles. 

\subsection{Adelic Divisors}\label{AdelicDivisors}

The adelic divisors and adelic line bundles on quasi-projective varieties are introduced by Yuan--Zhang in \cite{YZ}.  Roughly speaking, an adelic line bundle on a quasi-projective variety $ X $ is a reasonable limit of a sequence of (hermitian) line bundles on projective models of $ X $. They also introduces an intersection theory and Deligne pairing for adelic line bundles.

In \cite{YZ}, they mainly used adelic divisors on quasi-projective varieties over $ \mathbb{Z} $ or a field. They also defined adelic line bundles in a more general setting, including a more general valuation ring as a base, or Berkovich space. 

In our paper we care more on analytic spaces, and it is slightly different with the algebraic one. To explain the notations and emphasize the differences, we give a brief introduction in the analytic situations.

We will use both the algebraic notation and analytic notations. Throughout this paper, all schemes are locally of finite type over $ \C $. For a scheme or analytic space $ X $, let $ \Div(X) $ be the group of Cartier divisors on $ X $. Divisors will always mean Cartier divisors, unless stated otherwise. 

First we omit the hermitian metric. 

Let $ X $ be an reduced and irreducible scheme or an analytic space, and $ Z\subset X $ be a subvariety. Let $ U=X\backslash Z\subset X $ be the open subspace of $ X $. We define divisors of mixed coefficients
\begin{align*}
\Div(X,U) := \Div(X)_\Q \times_{\Div(U)_\Q} \Div(U).
\end{align*}
In other words, an element $ E\in \Div(X,U) $ is a pair $ E=(E_1,E_2) $ where $ E_1\in \Div(X)_\Q $ and $ E_2\in \Div(U) $ with the same image in $ \Div(U)_\Q $. By abuse of notations, we usually write $ E_1=E $ and $ E_2=E|_U $. 

The model divisors on $ (X,U) $ is defined to be
\begin{align*}
\Div(X,U)_{\mod} = \lim_{\substack{\longrightarrow\\X'}} \Div(X',U)
\end{align*}
where the limit is taken over all projective morphism $ \pi:X'\to X $, where $ X' $ is reduced and irreducible, and $ \pi^{-1}(U)\to U $ is an isomorphism. The morphisms are the pull-back of divisors. This definition depends on $ X $. Note that if $ X $ is itself projective over $ \C $, it coincides with the Yuan--Zhang's definition $ \Div(U/\C)_{\mod} $ in \cite[\S 2.4.1]{YZ}. 

If one wants to do completion as in \cite[\S 2.4.1]{YZ}, one realizes that in the analytic case, two different boundary divisors do not bound each other, as there are probably infinitely many irreducible components of an analytic divisor. Thus, to define adelic divisors in this case, we need to choose a fixed boundary divisor. 

By an analytic triple $ (X,U,Z)$, we mean an analytic variety $ X $, an analytic effective divisor $ Z\subset X $, and the open subspace $ U=X\backslash Z \subset X $. 

For each analytic triple $ (X,U,Z) $, we get a boundary norm on $ \Div(X,U)_{\mod} $ as
\begin{align*}
\|\cdot\|_{Z}: \Div(X,U)_{\mod} \to [0,\infty]
\end{align*}
by
\begin{align*}
\|{E}\|_{Z} = \inf\{ \varepsilon\in \Q_{>0}:-\varepsilon {Z} \leqslant {E}\leqslant \varepsilon {Z} \}
\end{align*}
where the partial order is given by effectivity. 

Define the analytic adelic divisor $ \Div(X,U,Z) $ to be the completion of $ \Div(X,U)_{\mod} $ with respect to the boundary norm $ \|\cdot\|_Z $. Again, if $ X $ is a projective variety, then $ Z $ is also algebraic by Serre's GAGA \cite{Ser}, and the completion procedure do not depend on $ Z $. In this case, our definition will be the same as Yuan--Zhang's $ \Div(U/\C) $. 

If $ X $ is normal, $ Z,Z' $ are two boundary divisors with finitely many irreducible components, and $ |Z|=|Z'| $, then there is an integer $ n>0 $ such that $ \frac{1}{n} Z\leqslant Z'\leqslant nZ $, hence they define the same completion and $ \Div(X,U,Z)\cong \Div(X,U,Z') $. This occurs especially when $ X $ is the unit disk $ \Delta $, and $ U=\Delta\backslash\{0\} $. 

Let $ (X,U,Z) $ be an analytic triple. Let $ Y $ be an analytic variety, and let $ f:Y\to X $ be a holomorphic map. Assume that $ \im(f)\not\subset |Z| $. Let $ W=f^*Z $ and $ V=Y\backslash W $. Then $ V\neq \emptyset $, and we get an analytic tuple $ (Y,V,W) $. 

By a morphism $ f:(X,U,Z)\to (Y,V,W) $ of analytic triples, we mean a holomorphic map $ f:Y\to X $, such that $ \im(f)\not\subset Z $, $ f^*Z=W $, and $ f^{-1}(U)=V $. If $ \dim \im(f)=\dim X $, then morphism $ f $ induces a pull-back map of adelic divisors
\begin{align*}
f^*:\Div(X,U,Z)\to \Div(Y,V,W).
\end{align*}
One example is that $ Y $ is an open subset of $ X $. Thus we have a well-defined restriction map to open subset. 

If the image does not have maximal dimension, the pull-back map does not exist for all divisors. We consider only the subset of $ \Div(X,U,Z) $ supported on the boundary. Define
\begin{align*}
	\Div(X,U,Z)_{b} = \{E\in \Div(X,U,Z):E|_U=0\}.
\end{align*}
Let $ f:(X,U,Z)\to (Y,V,W) $ be a morphism of analytic triples, then $ f $ induces a pull-back morphism
\begin{align*}
	f^*: \Div(X,U,Z)_b\to \Div(Y,V,W)_b.
\end{align*}

In the following, we mainly take $ Y $ to be the unit disk $ \Delta $, and $ f:\Delta\to X $ satisfies $ f(\Delta\backslash\{0\})\subset U $. Let $ W=f^*Z=c\cdot\{0\} $ for some $ c\geqslant 0 $, and $ V=f^{-1}(U) $. 

Then $  \Div(\Delta,V,W)_b  $ is of dimension $ 0 $ or $ 1 $, depending on whether $ c=0 $ or not. There is a canonical injection $ \Div(\Delta,V,W)_b \to \Div(\Delta, \Delta\backslash \{0\},\{0\})_b $, and we have a multiplicity map on $ \Div(\Delta,\Delta\backslash\{0\},\{0\})_b $ normalized so that $ \mult_0 \{0\}=1 $. By pulling back, $ f $ induces a group morphism
\begin{align*}
	\Div(X,U,Z)_b&\to \R\\
	E&\mapsto \mult_0 f^*E.
\end{align*}

Similar as in \cite[Def.5.1.1]{YZ}, an analytic adelic divisor is called effective if it can be represented by a Cauchy sequence of effective model divisors. We have the following proposition. 
\begin{prop}\label{EffectivityLemma}
	Let $ (X,U,Z) $ be an analytic triple with $ X $ smooth, and let $ E\in \Div(X,U,Z)_b $. Let $ \Delta $ be the unit disk. Then $ E $ is effective if and only if for any holomorphic map $ f:\Delta\to X $ with $ f(\Delta\backslash \{0\})\subset U $, we have 
	\begin{align*}
	\mult_0 f^* E\geqslant 0. 
	\end{align*}
\end{prop}

\begin{proof}
	If $ E $ is effective, then any pull-back of $ E $ is also effective, and thus has non-negative degree. So we only need to prove the converse. 
	
	Suppose that $ E $ is represented by a Cauchy sequence $ \{E_i\in \Div(X_i,U) \} $. We may assume that each $ X_i $ is smooth, $ E_i|_U=0 $, and there is a projective morphism $ X_i\to X_{i-1} $ compactible with the open immersions $ U\to X_i $. 
	
	Since $ X_i $ is smooth, $ E_i $ is an analytic Weil $ \Q $-divisor on $ X_i $. Let $ \{F_{ij}:j\in I_i\} $ be the set of prime boundary Weil divisors in $ X_i $, which are locally finite on $ X_i $. Then $ E_i=\sum_{j\in I_i}c_{ij} F_{ij} $ for suitable $ c_{ij}\in \Q $. 
	
	Now take any $ n>i $. Since $ X_{n}\to X_i $ is birational and projective and all $ X_i $ are smooth, the rational inverse map $ X_i\dashrightarrow X_{n} $ is defined outside a closed subvariety of codimension $ \geqslant 2 $. Thus there is a unique prime divisor $ F_{n,j_{n}} $ on $ X_{n} $ birational to $ F_{ij} $. Let $ c_{n,j_n} $ be the multiplicity of $ E_n $ at $ F_{n,j_n} $. 
	
	Since $ \{E_i\} $ is a Cauchy sequence, there exists $ \varepsilon_n\in \Q_{>0} $, such that $ \varepsilon_n\to 0 $, and
	\begin{align*}
	-\varepsilon_n Z\leqslant  E_m-E_n\leqslant \varepsilon_n Z\ \ \text{ for all }\ \  m>n.
	\end{align*} 
	Let $ d_{ij}\in \mathbb{Z}_{>0} $ be the multiplicity of $ Z $ at $ F_{ij} $. Compare the multiplicity at $ F_{n,j_n} $, we get
	\begin{align}\label{IneqEff}
		-\varepsilon_n d_{ij}\leqslant c_{m,j_m}-c_{n,j_n}\leqslant \varepsilon_n d_{ij}.
	\end{align}
	for all $ m>n\geqslant i $. Thus the limit $ m_{ij}=\displaystyle\lim_{n\to\infty}c_{n,j_n} $ exists. We call it the multiplicity of $ E $ at $ F_{ij} $. 
		
	Since there are only countably many varieties and divisors, and $ \C $ is uncountable, for a very general point $ x\in F_{ij} $, $ x $ lies on exactly one prime divisor on $ X_i $, and the inverse map $ X_i\dashrightarrow X_n $ is defined at $ x $ for all $ n>i $. 
	
	Given any morphism $ f:\Delta \to X_i $ such that $ f(0)=x\in F_{ij} $ and $ f(\Delta\backslash \{0\})\subset U $, we have
	\begin{align*}
	0&\leqslant \mult_0 f^*E=\lim_{n\to\infty} \mult_0 f^*E_{n} \\
	&= \lim_{n\to\infty} c_{n,j'}\cdot\mult_0 f^*F_{n,j'} =  m_{ij}\cdot\mult_0 f^*F_{ij}.
	\end{align*}
	Thus $ m_{ij}\geqslant 0 $.

	Now choose any $ m'_{ij}\in \Q_{\geqslant 0} $ such that $ |m'_{ij}-m_{ij}|\leqslant \varepsilon_i d_{ij} $. Define $ E'_i=\sum_{j\in I_j} m'_{ij} F_{ij}\in \Div(X_i,U) $. Then $ E_i' $ is effective.
	
	In the inequality (\ref{IneqEff}), let $ m\to \infty $ and $ n=i $, then 
	\begin{align*}
	-\varepsilon_i d_{ij}\leqslant m_{ij}-c_{n,j_n}\leqslant \varepsilon_i d_{ij}.
	\end{align*}
	Thus we have
	\begin{align*}
	-2 \varepsilon_i d_{ij}\leqslant m'_{ij}-c_{n,j_n}\leqslant 2\varepsilon_i d_{ij},
	\end{align*} 
	and thus $ -2\varepsilon_i Z\leqslant E_i'-E_i\leqslant 2\varepsilon_i Z $	
	Then $ \{E'_i\} $ is a Cauchy sequence of effective divisors converging to $ E $. Hence $ E $ is effective. 
\end{proof}

Now we consider the hermitian divisors. Let $ X $ be an algebraic or analytic variety. An hermitian divisor $ \bar{D}=(D,g_{\bar{D}}) $ on $ X $ consists a divisor $ D $ on $ X $ and a continuous Green function $ g_{\bar{D}}(x) $ on $ X\backslash |D| $, which induces a continuous hermitian metric $ \|\cdot\| $ on the line bundle $ \O(D) $ by
\begin{align*}
\|1\|(x)=e^{-g_{\bar{D}}(x)} ,
\end{align*} 
where $ 1 $ is the canonical section on $ \O(D) $. It is called effective if $ D $ is effective and $ g_{\bar{D}}\geqslant 0 $. We define $ \hDiv(X) $ to be the group of hermitian divisors on $ X $. 

Most definitions for hermitian adelic divisors are only slightly different from \cite{YZ} as above. 

Let $ X $ be an algebraic or analytic variety, and $ Z\subset X $ be a subvariety. The divisor of mixed coefficients is defined as

\begin{align*}
\hDiv(X,U) := \hDiv(X)_\Q \times_{\Div(U)_\Q} \Div(U),
\end{align*}
and the model hermitian divisor is defined as 
\begin{align*}
\hDiv(X,U)_{\mod} = \lim_{\substack{\longrightarrow\\X'}} \hDiv(X',U),
\end{align*}
where $ X' $ is taken over all projective morphism $ \pi:X'\to X $ such that $ X' $ is a variety and $ \pi^{-1}(U)\cong U $. 

A hermitian boundary divisor $ \bar{Z}=(Z,g_{\bar{Z}}) $ means an effective divisor $ Z $ and a strictly positive green function $ g_{\bar{Z}} $ of $ Z $. A hermitian analytic triple $ (X,U,\bar{Z}) $ means an analytic variety $ X $, an hermitian boundary divisor $ \bar{Z}=(Z,g_{\bar{Z}}) $, and the open subset $ U=X\backslash Z\subset X $. 

For each hermitian analytic triple $ (X,U,\bar{Z}) $, we det a boundary norm
\begin{align*}
	\|\cdot\|_{\bar{Z}}: \hDiv(X,U)_{\mod} \to [0,\infty]
\end{align*}
by
\begin{align*}
\|{\bar{E}}\|_{\bar{Z}} = \inf\{ \varepsilon\in \Q_{>0}:-\varepsilon \bar{Z} \leqslant \bar{E}\leqslant \varepsilon \bar{Z} \}.
\end{align*}
Let $ \hDiv(X,U,Z) $ be the completion of $ \hDiv(X,U)_{\mod} $ with respect to the boundary norm.

We emphasize that given two boundary divisors $ \bar{Z}_1,\bar{Z}_2 $, even if $ Z_1 $, $ Z_2 $ both have finitely many irreducible components, they do not bound each other in general. We only have the following results.
\begin{prop}\label{Bound}
	Let $ (X,U,\bar{Z}_1) $ and $ (X,U,\bar{Z}_2) $ are two analytic triples with $ |Z_1|=|Z_2| $. Assume that $ X $ is normal. Then for any open subset $ V\subset X $ such that the closure of $ V $ in $ X $ is compact, there is a positive integer $ n>0 $, such that
	\begin{align*}
	\frac{1}{n} \bar{Z}_1|_V\leqslant \bar{Z}_2|_V\leqslant n\bar{Z}_1|_V
	\end{align*}
\end{prop}
\begin{proof}
	First we choose an open neighborhood $ W $ of $ \bar{V} $, such that the closure $ \overline{W} $ in $ X $ is compact. 
	
	Since the closure of $ W $ is compact, there are only finitely many irreducible components on $ Z_i $ which intersects the closure of $ W $. Thus there exists a positive integer $ n'>0 $, such that 
	\begin{align*}
	\frac{1}{n'} {Z}_1|_W \leqslant {Z}_2|_W\leqslant n'Z_1|_W.
	\end{align*}

	Next we consider the Green functions. The formula above shows that 
	\begin{align}\label{4.5}
	\frac{1}{2n'} g_{\bar{Z}_1}(x) \leqslant g_{\bar{Z}_2}(x) \leqslant 2n' g_{\bar{Z}_1}(x)
	\end{align}
	for all points in $ x \in W\backslash |Z_1| $ that are sufficiently close to $ |Z_1| $. In particular, the formula (\ref{4.5}) holds for an open neighborhood $ V' $ of $ |Z_1|\cap \bar{V} $. 
	
	The rest part $ \bar{V}\backslash V' $ is a compact subset, and both Green functions are positive. Thus they bound each other by constants. 	
\end{proof}

\subsection{Admissible line bundles}
Let $ K=\C $ with Euclidean metric, or $ K $ be a complete field with a discrete valuation, and $ C $ be a smooth projective geometrically integral curve over $ K $. Let $ \omega_C $ be the dualizing sheaf on $ C $ and $ \Delta $ be the diagonal divisor on $ X\times X $. The admissible metrics on $ \omega_C $ and $ \O(\Delta) $ on $ C $ were introduced by Arakelov in \cite{Ara} when $ K $ is archimedean, and by Zhang in \cite{Zha93} when $ K $ is non-archimedean. We denote the metrics by $ \|\cdot\|_{a} $, $ \|\cdot\|_{\Delta,a} $, respectively. 

Yuan generalized it to family of smooth curves in \cite{Yua} using adelic line bundles. He used Berkovich space to state his result. For simplicity, we only state the archimedean case, i.e. $ K=\C $ with the usual Euclidean metric. 

\begin{theorem}[Yuan]	
	Let $ S $ be a quasi-projective normal integral scheme over $ \C $. Let $ \pi:X\to S $ be a smooth relative curve of genus $ g>0 $. Denote by $ \Delta:X\to X\times_S X $ be the diagonal morphism. 
	
	(1) There is an adelic line bundle $ \bar{\omega}_{X/S,a} $ in $ \widehat{\mathcal{P}ic}(X) $ with the underlying line bundle $ \omega_{X/S} $, such that for each $ s\in S $, the metric of $ \omega_{X_s} $ is equal to the canonical admissible metric $ \|\cdot\|_{a} $
	
	(2) There is an adelic line bundle $ \bar{\O}(\Delta)_{a} $ in $ \widehat{\mathcal{P}ic}(X\times_S X) $ with the underlying line bundle $ \O(\Delta) $, such that for each $ s\in S $, the metric of $ \O(\Delta)|_{X_s\times X_s} $ is equal to the canonical admissible metric $ \|\cdot\|_{\Delta,a} $. 
	
	Moreover, the adelic line bundles satisfy the following properties: 
	
	(1) The canonical isomorphism
	\begin{align*}
	\omega_{X/S}\to \Delta^* \O(-\Delta)
	\end{align*}
	induces an isometry
	\begin{align*}
	\bar{\omega}_{X/S,a}\to \Delta^* \bar{\O}(-\Delta)_a
	\end{align*}
	
	(2) The canonical isomorphism 
	\begin{align*}
	p_{1*}\langle \O(\Delta),p_2^* \omega_{X/S}\rangle \to \Delta^* p_2^* \omega_{X/S} \to \omega_{X/S} 
	\end{align*}
	induces an isometry
	\begin{align*}
	p_{1*}\langle \bar{\O}(\Delta)_a,p_2^* \bar{\omega}_{X/S,a} \rangle \to \bar{\omega}_{X/S,a}
	\end{align*}
	Hee $ p_1,p_2:X\times_S X\to X $ denote the two projections.  
\end{theorem}

\subsection{Zhang--Kawazumi's $ \phi $-invariants for Curves}\label{Zhang}
In \cite{Zha10}, Zhang introduced the $ \phi $-invariants when studying the height of Gross-Schoen cycles. He defined the $ \phi $-invariants of a smooth projective curve over both archimedean non-archimedean local fields. The archimedean $ \phi $-invariants was also independently introduced by Kawazumi in \cite{Kaw08,Kaw09}.

\subsubsection*{Zhang--Kawazumi's $ \phi $-invariants for Riemann Surfaces}

If $ K=\C $ and $ C $ is a compact Riemann surface of genus $ g\geqslant 2 $, the $ \phi $-invariants is defined as
\begin{align*}
\phi(C) = -\int_{(C\times C)} g_{\Ar}(x,y) c_1(\bar\O(\Delta)_a)^2.
\end{align*}

In \cite{Zha93}, it has explicit expressions as follows.
\begin{align*}
\phi(C)=\sum_{\lambda} \sum_{j=1}^g \sum_{k=1}^g \frac{2}{\lambda_n} 
\Big| \int_C \phi_\lambda \omega_j\wedge \bar{\omega}_k\Big|
\end{align*}
where $ \omega_1,\cdots,\omega_g $ is an orthonormal basis of $ \Gamma(C,\omega_C) $ with respect to the hermitian product
\begin{align*}
\langle \alpha,\beta \rangle =\frac{i}{2}\int_C \alpha\wedge\bar{\beta},
\end{align*}
and the first summation goes over all positive eigenvalues $ \lambda $ of the Laplacian operator
\begin{align*}
\Delta_{d\mu} f =\frac{1}{\pi i}\partial\bar{\partial}f / d\mu
\end{align*}
over $ C^{\infty} (C) $, and $ \phi_\lambda $ is an eigenvector of $ \lambda $ normalized such that $ \{\phi_\lambda\}_\lambda $ is orthonormal with respect to the inner product
\begin{align*}
\langle f_1,f_2 \rangle = \int_C f_1\bar{f}_2 d\mu.
\end{align*}

\subsubsection*{Zhang's $ \phi $-invariants for Non-archimedean Fields}

Now let $ K $ is a complete non-archimedean field with a discrete valuation, and $ \O_K $ be the ring of integers of $ K $, and $ k $ be its residue field. Let $ C $ be a smooth curve over $ K $ of genus $ g\geqslant 2 $. We first discuss the dual graph of $ C/K $.

By the semistable reduction theorem \cite[Cor. 2.7]{DM}, there exists a finite extension $ K'/K $, with ring of integers $ \O_{K'} $ and residue field $ k' $, such that $ C_{K'} $ has split semistable reduction, i.e., a semistable reduction satisfying that all nodes are defined over $ k' $ and all tangent vector at the nodes are defined over $ k' $. This automatically holds for any semistable reduction if $ k $ is algebraically closed. We will only use the case $ k=\C $ in this paper. 

Then we have the minimal regular model $ \mathcal{C}/\O_{K'} $, which is split semistable. The reduction graph of $ C/K $ is given as follows.

The underlying graph is just the dual graph of the special fiber of $ \mathcal{C}_k $, i.e., for each irreducible component $ D $ of $ \mathcal{C}_{k'} $, we have a vertex $ v(D) $; and for each node $ N $ in $ \mathcal{C}_{k'} $, we have an edge $ e(N) $ in $ \Gamma(C) $ connecting the vertices $ v(D_1),v(D_2) $, where $ D_1,D_2 $ are the irreducible components of $ \mathcal{C}_{k'} $ containing $ N $(possibly $ D_1=D_2 $). 

All edge have the same length $ \frac{1}{[K':K]} $, and the polarization $ q $ at the vertex $ v(D) $ is the geometric genus of the normalization of the irreducible component $ D $. 

Now the dual graph $ \Gamma(C) $ together with the metric and polarization, is a polarized metrized graph. The $ \phi $-invariant of $ C $ is just defined as the $ \phi $-invariant of its dual graph 
\begin{align*}
\phi(C)=\phi(\Gamma(C)).
\end{align*}

For later use, for a split semi-stable curve over the residue field $ k $, we also define its dual graph in the same way, and we assume that all edges have length $ 1 $.

\subsection{Globalization of $ \phi $-invariants}\label{GlobalizationofPhi}
Now we consider the globalization $ \bar{\Phi} $ of the $ \phi $-invariants defined in \cite[\S 3.3.2]{Yua}. Let $ (\pi,\pi):X\times_S X\to S $ be the structure morphism, and $ \Delta:X\to X\times_S X $ be the diagonal morphism. By the property of Deligne pairing, there are canonical isomorphisms
\begin{align*}
(\pi,\pi)_*\langle \O(\Delta),\O(\Delta),\O(\Delta)\rangle \to \pi_*\langle \O(\Delta),\O(\Delta)\rangle \to \pi_*\langle\omega,\omega \rangle
\end{align*}
which defines a section $ s $ of adelic line bundle
\begin{align*}
\pi_*\langle\bar{\omega}_{X/S,a},\bar{\omega}_{X/S,a} \rangle - (\pi,\pi)_*\langle \bar{\O}(\Delta)_a,\bar{\O}(\Delta)_a,\bar{\O}(\Delta)_a\rangle 
\end{align*}
The globalization of the $ \phi $-invariants is defined by
\begin{align*}
\bar{\Phi}=(\Phi,g_{\bar{\Phi}}) : = \widehat{\div}(s).
\end{align*}
Then we have
\begin{align*}
\O(\bar{\Phi}) = \pi_*\langle\bar{\omega}_{X/S,a},\bar{\omega}_{X/S,a} \rangle - (\pi,\pi)_*\langle \bar{\O}(\Delta)_a,\bar{\O}(\Delta)_a,\bar{\O}(\Delta)_a\rangle,
\end{align*}
and the underlying divisor $ \Phi|_S = 0 $. By the integration formula, we see that 
\begin{align*}
g_{\bar{\Phi}}(s) = \int_{X_s} \log \|1\|_{\Delta,a} c_1(\bar{O}(\Delta)_{a})^2 = \phi(X_s). 
\end{align*}

Finally, in \cite[\S 3.3.2]{Yua}, Yuan used the Berkovich space to state his result, and hence both archimedean and the non-archimedean $ \phi $-invariants appear in $ \bar{\Phi} $. We give an explanation. If $ S $ is a smooth curve over $ \C $, and $ \pi:X\to S $ admits a semistable compactification $ \bar{\pi}:\bar{X}\to \bar{S} $ and $ \bar{S} $ is a smooth projective curve. Assume that $ \bar{X} $ is smooth. Then the underlying divisor $ \Phi $ has the form
\begin{align*}
\Phi= \sum_{s\in \bar{S}\backslash S}\phi(\Gamma(X_s)) s,
\end{align*}
where $ X_s $ is the fiber of $ \bar{\pi} $ at $ s $. If $ \bar{X} $ is not smooth, it still holds except $ X_s $ should be replaced by $ \bar{X}\times \mathrm{Spec}\widehat{\O_{S,s}} $ as a smooth curve over $ \mathrm{Frac}\widehat{\O_{S,s}}$.

\section{Proof of Theorem \ref{PhiForRiemannSurface2}}\label{Proof2}
In this section we first review some basic properties of the Kuranishi family of stable curves in \cite[\S XI.4]{ACGH}, and then we prove Theorem \ref{PhiForRiemannSurface2}.

Let $ X_0 $ be a stable curve with nodes, and $ \pi:X\to S $ be a standard algebraic Kuranishi family of $ X_0 $, see \cite[Thm. XI.6.5]{ACGH}. In particular, $ S $ is an affine scheme, and $ \pi $ is Kuranishi at every point $ s\in S $. Then $ X $, $ S $ are both smooth by \cite[Thm. XI.3.17, Cor XI.4.6]{ACGH}. More precisely, for any $ s_0\in S $, and let $ p_1,\ldots,p_r $ are all nodes on the fiber $ X_{s_0} $, there exist local coordinates $ t_1,\cdots,t_n $ centered at $ s_0 $ on $ S $, and near the node $ p_j\in X_{s_0}\subset X $, the total space $ X $ locally has the form $ xy=t_j $. Hence $ X $ is smooth over $ \C $.

In the following, assume that $ X_0 $ is a stable curve with no non-separating nodes. By \cite[\S X.2]{ACGH}, $ \mathrm{Pic}^0(X_0) $ is an abelian variety.

Let $ \pi:X\to S $ be an algebraic Kuranishi family of $ X_0 $, and we assume that for any $ s\in S $, the fiber $ X_s $ has no non-separating nodes. Let $ \Pic_{X/S} $ be the relative Picard functor, and $ \Pic^0_{X/S} $ be the identity component of $ \Pic_{X/S} $. By \cite[Thm. 1 in \S 9.4]{BLR}, the functor $ \Pic^0_{X/S} $ is representable by a smooth separated scheme $ J $ which is semi-abelian over $ S $ with a canonical rigidified symmetric ample line bundle $ \O(\Theta) $. By the reason above, it is in fact abelian over $ S $. We denote it by $ J/S $. Let $ e:S\to J $ be the identity section. 

Since $ \Theta $ is symmetric, $ [2]^* \Theta = 4\Theta $. The rigidification $ e^*\O(\Theta)\cong  \O(S) $ induces an isomorphism $ f:[2]^* \Theta \to 4\Theta $. Using Tate's limiting argument, there is a unique hermitian adelic line bundle $ \bar{\Theta} = (\tilde{\Theta},\|\cdot\|_{\Theta}) $ extending $ \Theta $, such that $ f $ becomes an isometry. See \cite[\S 6.1]{YZ} for details. 

Let $ U\subset S $ be the open subset parametrizing smooth curves, $ Y=\pi^{-1}(U) $. By the universal property of relative Picard functor, there is a morphism 
$$ i_\omega:Y\to J, x\mapsto (2g-2)x-\omega_{Y/U}. $$
Note that if we work on $ X/S $, the divisor $ (2g-2)x-\omega_{X/S} $ is not algebraically equivalent to $ 0 $, so this formula only defines a rational map $ X\dashrightarrow J $. But since $ X $ is smooth, by \cite[\S 8.4, Cor. 6]{BLR}, the rational morphism $ i_\omega $ extends to a morphism $ X\to J $, which we still denote it by $ i_\omega $. 

By \cite[Thm. 2.10]{Yua}, we have the following formulas.
\begin{align*}
\pi_{*}\langle i_{\omega}^* \bar{\Theta},i_{\omega}^* \bar{\Theta} \rangle = 16g(g-1)^3 \pi_{*}\langle \bar{\omega}_{Y/U,a},\bar{\omega}_{Y/U,a} \rangle\ \  \text{ in}\ \hPic(U)_\Q
\end{align*}
\begin{align*}
i_{\omega}^* \bar{\Theta} = 4g(g-1)\bar{\omega}_{Y/U,a} - \pi_{*}\langle \bar{\omega}_{Y/U,a},\bar{\omega}_{Y/U,a} \rangle\ \ \text{in}\ \hPic(Y)_\Q.
\end{align*}

Summarizing, we get
\begin{align*}
\bar{\omega}_{Y/U,a}=\frac{1}{4g(g-1)} \left( i_\omega^*\bar{\Theta} -\frac{1}{16g(g-1)^3} \pi_{2,*} \langle i_\omega^*\bar{\Theta},i_\omega^*\bar{\Theta} \rangle\right)
\end{align*}
in $ \hPic(Y)_\Q $. However, the right hand side of above formula defines a hermitian $ \Q $-line bundle on $ X $. Hence we've showed the following result. 
\begin{prop}
	There is a hermitian $ \Q $-line bundle on $ X $, whose restriction to $ Y $ equals to $ \bar{\omega}_{X/S,a} $ as hermitian $ \Q $-line bundles.
\end{prop}

Next we deal with $ \bar{\O}(\Delta)_a $. Let $ j:Y\times_U Y\to J $, $ (x,y)\mapsto (x-y) $, and let $ \tilde{j}=[2g-2]\circ j $. View $ \tilde{j} $ as a rational map $ X\times_S X\dashrightarrow J $. Note that $ \tilde{j}= i_{\omega}\circ p_{1}- i_{\omega}\circ p_{1} $ on $ Y\times_U Y $, thus $ \tilde{j} $ extend to a morphism $ \tilde{j}:X\times_S X\to J $. 

By \cite[Thm. 2.10]{Yua}, we have the following formula.
\begin{align*}
j^*\bar{\Theta} = 2\bar{\O}(\Delta)_a + p_1^* \omega_{X_0/S_0,a}+ p_2^* \omega_{X_0/S_0,a}\ \ \text{in}\ \hPic(Y\times_U Y)_\Q.
\end{align*}
Note that $ [n]^* \bar{\Theta}= n^{2}\bar{\Theta} $, hence we have the formula
\begin{align*}
\bar{\O}(\Delta)_a= \frac{1}{2}\left(\frac{1}{(2g-2)^{2}} \tilde{j} ^* \bar{\Theta} -p_1^*\bar{\omega}_a -p_2^*\bar{\omega}_a \right)\ \ \text{in}\ \hPic(Y\times_U Y)_\Q.
\end{align*}
Hence we get the following proposition. 
\begin{prop}
	There is a hermitian $ \Q $-line bundle on $ X\times_S X $, whose restriction to $ Y\times_U Y $ equals to $ \bar{\O}(\Delta)_a $ as hermitian $ \Q $-line bundles.
\end{prop}

The globalization of the $ \phi $-invariants is defined as
\begin{align*}
\O(\bar{\Phi}) =\pi_{*}\langle \omega_{X/S,a},\omega_{X/S,a} \rangle - (\pi,\pi)_*\langle\bar\O(\Delta)_a,\bar\O(\Delta)_a, \bar\O(\Delta)_a\rangle . 
\end{align*}
\begin{corollary}
	There is a hermitian $ \Q $-line bundle on $ S $, whose restriction to $ U $ equals to $ \O(\bar{\Phi}) $ as hermitian $ \Q $-line bundles.
\end{corollary}
In other words, there exists $ n>0 $, such that $ n\bar{\Phi} = (n\Phi, ng_{\bar{\Phi}}) $ extends to an hermitian divisor on $ S $. Since $ \Phi|_U=0 $, we see that $ \Phi|_S $ is a rational combination of boundary divisors. 

Locally near $ s_0 $, if we choose a suitable coordinate, the boundary divisor is defined by the equation $ t_1\cdots t_r $, where $ r $ is the number of nodes in $ X_0 $. Let $ (g_i,g-g_i) $ be the type of the $ i $-th nodes. Note that given any curve in $ S $ intersecting the singular locus properly, the underlining divisor $ \Phi $ is given by the non-archimedean $ \phi $-invariant, which is the $ \phi $-invariant of the dual graph. By \cite[4.4.1]{Zha10}, if $ \Gamma $ is a tree, the $ \phi $-invariant is given by
\begin{align*}
\phi(\Gamma)=\sum_{j=1}^{[g/2]} \frac{2j(g-j)}{g} \delta_j,
\end{align*}
where $ \delta_i $ is the length of edges in $ \Gamma $ of type $ j $, which in our case is the number of nodes of type $ j $. Thus, we get
\begin{align*}
\Phi|_S = \sum_{j=1}^{r} \frac{2j(g-j)}{g} \Delta_j
\end{align*}
where $ \Delta_j $ is the boundary divisor on $ S $ of type $ i $, see \cite[p.339]{ACGH}

Finally, note that $ g_{\bar{\Phi}} $ is a Green function of the $ \Q $-divisor $ \Phi|_S $ in the usual sense, we get Theorem \ref{PhiForRiemannSurface2} for the Kuranishi family.  

Note that any family is locally pull-back of the Kuranishi family, we get Theorem \ref{PhiForRiemannSurface2} for all family $ \pi:X\to S $.

\section{Special Adelic Divisors}\label{Special}
In this section, we give a condition of adelic divisors such that their behaviors are easily to describe. We use the analytic condition in this section. 

Let $ S= \{( t_1,\ldots,t_n)\in \C^n:|t_i|<1\} $ be the polydisk. Let $ D_i $ be the Cartier divisor defined by $ t_i=0 $. Fix $ 1\leqslant r\leqslant n $, let $ D=\sum_{1\leqslant i\leqslant r} D_i $ and $ U=S\backslash D $. 

Recall that $ \Div(S,U,D)_b $ is the set of all adelic divisors $ E\in \Div(S,U,D) $ with $ E|_U=0 $. 

\begin{definition}
	An adelic divisor $ E \in \Div(S,U,D)_b $ is called special if there exists a continuous function $ f:\R_{\geqslant 0}^{r} \to  \R  $, homogeneous of degree $ 1 $, satisfying the following.
	
	Let $ \Delta $ be the unit disk, and $ p:\Delta\to S $ be any holomorphic morphism with $ p(\Delta\backslash \{0\})\subset U $. Let $ m_i=\mult_0 p^* D_i\in \Z_{\geqslant 0} $. Then we have
	\begin{align*}
	\mult_0 p^* E = f(m_1,\ldots,m_r).
	\end{align*}

\end{definition}
Here a function is homogeneous of degree $ 1 $ if for all real numbers $ m_1,\ldots,m_r\geqslant 0 $ and $ \lambda\geqslant 0 $, 
\begin{align*}
f(\lambda m_1,\ldots,m_r)=\lambda f(m_1,\ldots,m_r),
\end{align*}
and the space of all special divisor is denoted by $ \Div(S,U,D)_{s} $. 

Such a function $ f $ is called the \emph{skeletal function} of the adelic divisor $ E $. 

Note that since $ f $ is homogeneous of degree $ 1 $, it is determined by its resitriction on the standard $ (r-1) $-simplex
\begin{align*}
\mathcal{B}=\{ (m_1,\ldots,m_r)\in \R^r:m_i\geqslant 0,\sum_{i=1}^r m_i =1 \}.
\end{align*}
By abuse of language, the function $ f|_\mathcal{B} $ is also called the skeletal function of $ E $. 
Our main theorem in this section is

\begin{theorem}\label{5.2}
	Given any continuous function  $ f:\mathcal{B}\to \R $, there exists a unique special adelic divisor $ E\in \Div(S,U,D)_s $, such that $ f $ is the skeletal function of $ E $.

\end{theorem}

\subsection{Model Case}\label{ModelCase}
We sill construct many special model divisors so that they are dense in $ \Div(S,U,D)_s $. FIrst we give some notations of multi-indices. 

For each multi-index $ \mathbb{\vec{\alpha}}=(\alpha_1,\ldots,\alpha_r)\in \Z_{\geqslant 0}^r $, let $ t^{\vec{\alpha}} $ be the monomial $ t_1^{\alpha_1}\cdots t_r^{\alpha_r} $. For $ \vec{m}=(m_1,\cdots,m_r)\in \R^r $, let $ \vec{\alpha}\cdot  \vec{m}=\sum_{i=1}^r \alpha_i m_i $.

Let $ T\subset \Z_{\geqslant 0}^r\backslash\{\vec{0} \} $ be a finite subset of multi-indices, and let $ \mathcal{I}\subset \O_S $ be the coherent ideal sheaf generated by monomials $ t^{\vec{\alpha}} $ for all $ \alpha\in T $. Let $ Z=V(\mathcal{I}) $ be the analytic subspace defined by $ \mathcal{I} $ which is possibly non-reduced. Then $ |Z|\subset |D| $. 

Let $ \pi:X_T\to S $ be the blowing up of $ S $ along the subscheme $ Z $, and let $ E_T $ be the exceptional divisor in $ X_T $, which is the pull-back of $ Z $ in $ X_T $. Let $ \Delta $ be the unit disk. Then each holomorphic map $ p:\Delta\to S $ with $ p^{-1}(U)=\Delta\backslash \{0\} $ lifts to $ p':\Delta\to X_T $ by the properness of the blowing up. 

Let $ m_i=\mult_0 p^*D_i\in \Z_{\geqslant 0} $, and let $ \vec{m} = (m_1,\ldots,m_r) $. Then we have
\begin{align*}
\mult_0 p'^*E_T=\mult_0 p^* Z=\min_{\vec{\alpha} \in T}\{\vec{\alpha} \cdot \vec{m}\}.
\end{align*}

The functions $ \min_{\vec{\alpha}\in T} \{\vec{\alpha}\cdot \vec{m}\} $, viewed as functions of $ \vec{m}\in \R_{\geqslant 0}^r $, are continuous and homogeneous of degree $ 1 $. Thus $ X_T $ is a special model adelic divisors. Define $ \Div(S,U,D)_{s,\mod}' $ be the set of such adelic divisors. Note that $ \Div(S,U,D)'_{s,\mod} $ is just a set, not a group.

Then any $ \Q $-linear combination of such divisors is a model special divisor. Although it is not known whether any model special divisor can be get this way, we still write $ \Div(S,U,D)_{s,\mod} $ to denote the subspace of $ \Q $-linear combination of divisors in $ \Div(S,U,D)'_{s,\mod} $.

For each $ E\in \Div(S,U,D)_{s,\mod} $, let $ f_E $ be the skeletal function of $ E $. Let $ V' $ be the set of all skeletal functions of divisors in $ \Div(S,U,D)'_{s,\mod} $, and let $ V $ be the set of all skeletal functions of divisors in $ \Div(S,U,D)_{s,\mod} $. 

\begin{prop}
	For each $ f_1,f_2\in V $, we have:
	
	(1) For each rational numbers $ a,b\in \Q $,  $ af_1+bf_2\in V $. 
	
	(2) The minimum function $ \min\{f,f_2\} $ is in $ V $. 
\end{prop}
\begin{proof}
	Part (1) is clear, so it suffices to prove (2). We first assume $ f_1,f_2\in V' $. 
	
	Let $ Z_1,Z_2 $ be two closed subspaces of $ S $ defined by monomials, let $ X_i $ be the blow up of $ Z_i $ with exceptional divisor $ E_i $, and $ f_i $ is the skeletal function of $ E_i $. 
	
	Let $ X $ be the blow up of $ X_1 $ for the subspace $ \pi_1^* Z_2 $. Then the morphism $ X\to S $ factors through $ X_2 $ by the universal property of blowing ups. Now $ E_i $ is viewed as a divisor of $ X $ via pull-back. 
	
	By direct computation, the skeletal function of $ Z_1\cap Z_2 $ is $ \min\{f_1,f_2\} $, and the skeletal function of the closed subspace defined by the ideal sheaf $ \mathcal{I}_1\cdot \mathcal{I}_2 $ is $ f_1+f_2 $. Thus $ V' $ is closed under addition and taking minimum. 
	
	Now we consider taking minimum in $ V $. Given any $ f,g\in V $, suppose that $ f=f_1-f_2 $, $ g=g_1-g_2 $, where $ f_i,g_i $ are all $ \Q_{>0} $-linear combination of divisors in $ V' $. Then we have.
	\begin{align*}
	\min\{f,g\}=\min \{f_1-f_2,g_1-g_2\}=\min \{ f_1+g_2,g_1+f_2 \} -f_2-g_2\in V.
	\end{align*}
	Thus $ V $ is closed under taking minimum. 
\end{proof}

\subsection{Taking Limit}\label{TakingLimit}

In this section, We will show that any special adelic divisor is the limit of a sequence of divisors in $ \Div(S,U,D)_{s,mod} $. 

First, we give a variant of Stone-Weierstrass theorem. The proof is exactly the same as the proof of \cite[Thm. 7.32]{Rud}. 

\begin{prop}
	Let $ X $ be a compact Hausdorff space, and let $ C(X) $ be the Banach space of real-valued continuous functions on $ X $ with the uniform norm. Let $ A\subset C(X) $ be a $ \Q $-linear subspace of $ C(X) $ satisfying the followings:
	
	(1) The constant function $ 1\in A $.
	
	(2) The space $ A $ separates points, i.e., for each $ x,y\in X $, there exists a function $ f\in A $, with $ f(x)\neq f(y) $. 
	
	(3) For each pair $ f,g \in A $, we have $ \min\{f,g\}\in A $.
	
	Then $ A $ is dense in $ C(X) $. 
\end{prop}
Thus we get the following corollary.
\begin{corollary}
	The skeletal functions of divisors in $ \Div(S,U,D)_{s,\mod} $ is dense in $ C(\mathcal{B}) $. 
\end{corollary}

Now we consider taking limits for both adelic divisors and continuous functions.

\begin{prop}
	Let $ E_i\in \Div(S,U,D)_{s,\mod} $ be a sequence of special model divisors. Then $ E_i $ is a Cauchy sequence in the sense of adelic divisors if and only if the corresponding sequence of skeletal functions $ f_i $ is a Cauchy sequence in $ C(\mathcal{B}) $.
\end{prop}

\begin{proof}
	Recall that $ D=V(t_1\cdots t_r) $ is the boundary divisor we choose. Then $ E_i $ is a Cauchy sequence in the sence of adelic divisors if and only if there exist a sequence $ \varepsilon_i>0 $ with $ \varepsilon_i\to 0 $, such that 
	$$ -\varepsilon_i D\leqslant  E_j-E_i\leqslant \varepsilon_i D $$
	for any $ j\geqslant i $. 
	
	By Proposition \ref{EffectivityLemma}, the inequality above transforms as 	
	\begin{align}\label{EqEff}
	-\varepsilon_i \mult_0 p^*D\leqslant  \mult_0 p^*(E_j-E_i)\leqslant \varepsilon_i \mult_0 p^*D 
	\end{align}	
	for any $ j\geqslant i $ and all holomorphic map $ p:\Delta\to S $ with $ p(\Delta\backslash \{0\})\subset U $. 
	
	Let $ M=\mult_0 p^*D $. If $ M=0 $, then $ p(0)\in U $, and all terms in (\ref{EqEff}) are zero. Now assume $ M>0 $. Let $ m_i=\mult_0 p^* D_i $. Then the inequality (\ref{EqEff}) becomes
	\begin{align*}
	-\varepsilon_iM \leqslant M\left(f_j\left(\frac{m_1}{M},\ldots,\frac{m_r}{M}\right)-f_i\left(\frac{m_1}{M},\ldots,\frac{m_r}{M}\right) \right)\leqslant \varepsilon_i M,\
	\end{align*}
	for all $j\geqslant i$, which is equivalent to the condition that $ f_j(x) $ is a Cauchy sequence in $ \mathcal{B} $ as rational points are dense in $ \mathcal{B} $. 
\end{proof}

Now we can prove Theorem \ref{5.2}.
\begin{proof}
	Given any continuous function $ f $, there exist a sequence of  $ f_i\in V $ such that $ f_i $ converges to $ f $ uniformly. Then the corresponding divisors $ E_i $ is a Cauchy sequence in the sense of adelic divisors. Thus $ E=\lim E_i $ exists, and the skeletal function of $ E=\lim E_i $ is $ f $. 
	
	The uniqueness of the adelic divisor comes from Proposition \ref{EffectivityLemma}. 
\end{proof}

This construction of adelic divisors has an advantage that we can write a Green function explicitly. 

Note that the function 
$$ g_{\bar{D}}(t_1,\ldots,t_n)=-\sum_{1\leqslant i\leqslant r}\log |t_i|  $$
is a positive Green function of $ D $ on $ S $, and hence $ \bar{D}=(D,g_{\bar{D}}) $ is a hermitian boundary divisor on $ S $.

\begin{theorem}\label{Theorem5.7}
	If $ E $ is a special adelic divisor with skeletal function $ f:\mathcal{B}\to \R $. Extend $ f $ to a function on $ \R_{\geqslant 0}^r $ homogeneously of degree $ 1 $. Let 
	\begin{align*}
	g_{E}(t_1,\ldots,t_n) = f(-\log |t_1|,\ldots,-\log |t_r|).
	\end{align*}
	Then we have $ \bar{E}=(E,g_E)\in \hDiv(S,U,\bar{D}) $. 
\end{theorem}

\begin{proof}
	First, we note that if $ \mathcal{I}=(u_1,\ldots,u_m) $ is an ideal sheaf on $ S $, and $ X $ is the blow up of $ S $ along $ \mathcal{I} $ with exception divisor $ E_{\mathcal{I}} $, then the function 
	\begin{align*}
	g_{E_{\mathcal{I}}}(t_1,\ldots,t_n) = \min \{-\log |u_1|,\ldots,-\log |u_m| \}
	\end{align*}
	is a Green function of $ E_{\mathcal{I}} $ on $ X $. 
	
	If $ E\in \Div(S,U,D)_{s,\mod} $ , then we define the function
	\begin{align*}
	g_E(t_1,\ldots,t_n) = f_E( -\log |t_1|,\ldots,-\log |t_r| ).
	\end{align*}
	If $ E\in \Div(S,U,D)'_{s,\mod} $ is defined by the ideal sheaf $ t^{\vec{\alpha}} $ for $ \vec{\alpha} $ in the finite set $ T\subset \Z^r_{\geqslant 0}\backslash\{\vec{0}\} $, then we have
	\begin{align*}
	f_E(m_1,\ldots,m_r)=\min_{\vec{\alpha}\in T}\left\{\sum_{i=1}^r \alpha_i m_i \right\}
	\end{align*}
	 and we see that $ g_E $ is a Green function of $ E $. 
	
	We will show that this particular Green function converges to our form.
	
	Let $ E_i\in \Div(S,U,D)_{s,\mod} $ converges to $ E $. Then the corresponding skeletal functions $ f_i $ converge to $ f $  uniformly on $ \mathcal{B}$. Recall that $ |t_i|<1 $ and our choice of green function is  
	$$ g_{\bar{D}}(t_1,\ldots,t_n) =-\sum_{1\leqslant i\leqslant r} \log |t_i|. $$
	Then we have	
	\begin{align*}
	&\ \left|\frac{f_i(-\log |t_1|,\ldots,-\log |t_r|) - f(-\log |t_1|,\ldots,-\log |t_r|)}{g_{\bar{D}}(t_1,\ldots,t_n)}\right|\\
	&=  \left |f_i\left(\frac{-\log |t_1|}{g_{\bar{D}}(t_1,\ldots,t_n)},\ldots, \frac{-\log |t_r|}{g_{\bar{D}}(t_1,\ldots,t_n)}\right) - f\left(\frac{-\log |t_1|}{g_{\bar{D}}(t_1,\ldots,t_n)},\ldots, \frac{-\log |t_r|}{g_{\bar{D}}(t_1,\ldots,t_n)}\right)\right|.
	\end{align*}
	Thus $ (E_i,g_{E_i}) $ converges to $ (E,g_E) $ in $ \hDiv(S,U,\bar{D}) $, and hence $ g_E $ is a Green function of $ E $. 
\end{proof}

\section{Specialness of Admissible Adelic Divisors}
In this section we show that many adelic divisors discussed in the \ref{AdmissibleLineBundles} are special in suitable sense.

\subsection{Specialness of $ \Phi $ and Proof of the Theorem \ref{PhiForRiemannSurface1}}\label{Proof1}
We first show that the adelic divisor $ \Phi $ on the moduli space is special, and then prove Theorem \ref{PhiForRiemannSurface1}. 

Since any family of stable curves is locally pull back of the Kuranishi family, we may only consider the Kuranishi family in this section.

Let $ X_0 $ be any stable curve with $ r $ nodes $ x_1,\ldots,x_r $, and let $ \pi:X\to S $ be the algebraic Kuranishi family of $ X_0 $, with a distinguished point $ s_0\in S $ parametrizing $ X_0 $ . Let $ U\subset S $ parametrizing the smooth curves. Then $ \bar{\Phi}=(\Phi,g_{\bar{\Phi}}) $ is the adelic divisor on $ U $. We fix a boundary divisor $ (\bar{S}, Z,g_{\bar{Z}}) $, where $ \bar{S} $ is a  projective model of $ S $, and $ Z $ is an effective divisor on $ \bar{S} $ satisfying $ |Z|=\bar{S}\backslash U $, and $ g_{\bar{Z}} $ is a positive Green function of $ Z $ on $ U $. Let $ \bar{Z}= (Z,g_{\bar{Z}}) $. 

Now we view $ (\bar{S},U,\bar{Z}) $ be an analytic triple. Then the adelic divisor $ \bar{\Phi} $ is an analytic adelic divisor. Since $ S $ is open in $ \bar{S} $, $ \bar{\Phi} $ restricts to an analytic adelic divisor on $ (S,U,\bar{Z}|_S) $. 

In the analytic setting, we may shrink $ S $ such that $ S $ is a polydisk of dimension $ n=3g-3 $, centered at $ s_0 $ with local coordinate $ \{(t_1,\ldots,t_n):|t_i|< 1\} $, such that $ t_i=0 $ parametrizes the deformation which are locally trivial at $ x_i $. Let $ D=V(t_1\cdots t_r) $, and let $ D_i=V(t_i) $. Let $ U=S\backslash D $. Then $ D $ is the singular locus of $ \pi $. It's easy to see that $ Z|_S $ and $ D $ induces the same boundary topology on $ S $, as they both are divisors supported in $ |D| $, and $ |D| $ has only finitely many irreducible components. Thus $ \Phi\in \Div(S,U,D) $. 

Let $ \Delta $ be the unit disk, and let $ p:\Delta\to S $ be any holomorphic map with $ p^{-1}(U)=\Delta\backslash \{0\} $, and let $ s=p(0) $. The fiber product $ Y=X\times_S \Delta $ is a stable curve over $ \Delta $ with central fiber $ Y_0=X_s $. The dual graph of $ Y/\Delta $ is a polarized metrized graph.

Let $ (\Gamma,q) $ be the polarized graph of the fiber $ X_0 $. Let $ m_i=\mult_0 p^* D_i\in \mathbb{Z}_{\geqslant 0} $. The relation between dual graph of $ Y/\Delta $ and $ X_0 $ is as follows. 
\begin{prop}
	The dual graph of $ Y/\Delta $ is exactly the polarized metrized graph $ (\Gamma,q;m_1,\ldots,m_n) $	defined in section \ref{DegenerationofGraphs}.
\end{prop}
\begin{proof}
First we omit the metric, and consider the relation of the two fibers $ X_0 $ and $ Y_0=X_s $. Suppose $ (a_1,\ldots,a_n) $ is the coordinate of $ s\in S $. Let $ I=\{i:1\leqslant i\leqslant r, a_i\neq 0\} $, and let $ E\subset \Gamma(X_0) $ be the subgraph containing edges $ \{e_i\}_{i\in I} $. 

For $ 1\leqslant i\leqslant r $, the coordinate $ t_i $ controls the local desingularization at the node $ p_i $. If $ t_i\neq 0 $, the node disappear, and the two irreducible components containing the node $ p_i $(possibly equal) become one irreducible component. Then we have
\begin{align*}
\Gamma(X_s) = \Gamma(X_{s_0})/E. 
\end{align*}

Next we consider the polarization of the graph. Suppose that $ F' $ is an irreducible component of $ X_s $, and let $ \xi' $ be the corresponding vertex in the dual graph $ \Gamma(X_s) $, then we have
\begin{align*}
q'(\xi') = g(F'),
\end{align*}
where $ g(F') $ is the genus of the normalization of $ F' $.

To compute $ g(F') $, let 
$$ S'=\{(t_1,\ldots,t_n): t_i=0 \text{ for all }i\notin I, i\leqslant r,\ \text{and }t_i\neq 0 \text{ for all }i\in I \} $$
be the locally closed subspace of $ S $ containing $ s $, and let $ X'=X\times_S S' $ be the fiber product. Then $ X'/S' $ is a topological fiber bundle, and thus for any point in $ S' $, the fiber has the same reduction graph, and the deformation does not change any nodes. Thus we may take an irreducible component $ W' $ of $ X' $ whose restriction to $ X_s $ is $ F' $. 

Let $ W $ be the closure of $ W' $ in $ X $, then $ W\to \bar{S}' $ is a family of nodal curve, and hence the arithmetic genus of the fibers is a constant. Let $ F_1,\ldots,F_m $ be the irreducible components on $ Y_{s_0} $, then we get
\begin{align*}
p_a (F') = p_a(\cup_{i=1}^m F_j).
\end{align*}

Now we can compute $ g(F') $ as follows. 
\begin{align*}
g(F') & = p_a(F')-\text{ number of nodes in } F\\
&= p_a(\cup_{i=1}^m F_j) -\text{ number of nodes in } F\\
& = \sum_{i=1}^m g(F_i)+1-m+\text{ number of nodes in }\cup F_i\\
&= \sum_{i=1}^m g(F_i) +1-m+|I|\\
&=  \sum_{i=1}^m g(F_i)+g(E_0)
\end{align*}
where $ E_0\subset \Gamma(X_{s_0}) $ is the preimage of $ \xi'\in \Gamma(X_s) $ under the quotient map $ \Gamma(X_{s_0})\to \Gamma(X_s) $. Thus the polarization on $ \Gamma/E $ coincide with the definition \ref{Def2.4}. 

Finally, the length of each edge is just the intersection number by a standard desingulariation of $ A_n $ surface singularities. 
\end{proof}

\begin{corollary}
	Let $ (\Gamma,q) $ be the reduction graph of $ X_0 $, and let 
	$$ \phi(\Gamma,q;L_1,\ldots,L_r) $$
	be the $ \phi $-function of graphs. Then for any holomorphic map $	p:\Delta\to S $ with $ p(\Delta\backslash \{0\})\subset U $, let $ m_i=\mult_0 p^*D_i $, then we have
	\begin{align*}
	\mult_0 p^*\Phi = \phi(\Gamma,q;m_1,\ldots,m_r)
	\end{align*}
	Thus $ \Phi $ is a special divisor with skeletal function $ \phi(\Gamma,q;L_1,\ldots,L_r) $. 
\end{corollary}

Now we can prove our main Theorem \ref{PhiForRiemannSurface1}.
\begin{proof}
	By Theorem \ref{Theorem5.7}, the function 
	$$ g'_{\Phi}(t_1,\ldots,t_r) = \phi(\Gamma(X_0),q;-\log |t_1|,\ldots,-\log |t_r|) $$
	is a Green function of $ \Phi\in \Div(S,U,D) $ with respect to the boundary divisor $ (D,g_D=-\sum_{i=1}^r \log |t_i| ) $. 
	
	Note that $ \bar{\Phi}\in \hDiv(S,U,\bar{Z}) $. Although $ g_Z $ and $ g_D $ may not bound each other on the whole $ S $ by a constant, they do bound each other on $ S'=\{(t_1,\ldots,t_n):|t_i|\leqslant \delta\} $ for a real number $ \delta<1 $ by Proposition \ref{Bound}. Let $ U'=S'\backslash |D| $. Then $ \bar{\Phi}|_{S'}\in \hDiv(S',U',\bar{D}) $. 
	
	Now $ (\Phi,g_{\bar{\Phi}}) $ and $ (\Phi,g'_{\Phi}) $ are both in $ \hDiv(S',U',\bar{D}) $. By \cite[Thm. 3.6.4]{YZ}, any two Green functions with the same underlying divisors differ at most $ o(g_D) $, thus we get our main theorem. 
\end{proof}

\subsection{Other Specialness}\label{OtherResults}
In this section we discuss the specialness of adelic divisors induced from the admissible line bundles $ \omega_a $ and $ \O(\Delta)_a $. 

Let $ X $ be a smooth variety, $ Z\subset X $ be a divisor of normal crossing, $ U=X\backslash Z $, and $ D\in\Div(X,U,Z)_b $. We say that $ D $ is a special divisor if it is special at all points $ x\in X $. 

Let $ \pi:X\to S $ be the Kuranishi family of stable curves as before. Let $ U\subset S $ parameterizing smooth curves, $ Y=\pi^{-1}(U) $, $ D\subset S $ be the boundary divisor, and let $ Z=\pi^*D $. Then $ Z $ is a divisor of normal crossing.

Consider two adelic line bundles $ \omega_{X/S},\omega_{Y/U,a}\in \Pic(X,Y,Z) $. The globalization of Zhang's $ \varepsilon $-invariant is defined in \cite[\S 3.2]{Yua} as 
\begin{align*}
\O(E)=\pi_*\langle \omega_{X/S},\omega_{X/S}\rangle-\pi_*\langle \omega_{Y/U,a},\omega_{Y/U,a}\rangle.
\end{align*}
We explain the formula as follows. The underlying line bundle of the right-hand side over $ U $ is canonically isomorphic to the trivial line bundle $ \O_U $. The section $ 1\in \O_U $ gives a rational section of the right-hand side, and the adelic divisor $ E:=\div(1) $ has underlying divisor $ 0 $ on $ S $. Similar as the $ \phi $-invariant, we get
\begin{prop}
	The adelic divisor $ E $ is special. 
\end{prop}

Similar as above, the identity morphism of $ \omega_{Y/U} $ induces as adelic divisor $ G\in \Div(X,Y,Z) $ with underlying divisor $ 0 $ on $ Y $ such that 
\begin{align*}
\O(G) = \omega_{X/S} - \omega_{Y/U,a}.
\end{align*}

Next we consider the $ \O(\Delta)_a $. The product $ X\times_S X $ is singular at the product of nodes. Let $ T $ be blow up $ X\times_S X $ at the singular locus, and $ (\pi,\pi):T\to S $ be the structure morphism. Then $ T $ is smooth, and $ (\pi,\pi)^* D $ is a divisor of normal crossing. 

Let $ \tilde{\Delta} $ be Zariski closure of the diagonal divisor $ \Delta\in \Div(Y\times_U Y) $ in $ T $. Then $ \O(\tilde{\Delta}) $ and $ \O(\Delta)_a $ both lie in $ \Div(T,Y\times_U Y,(\pi,\pi)^*D) $ with the same underlying divisor $ \Delta\in \Div(Y\times_U Y) $. The identity morphism of $ \O(\Delta)\in \Div(Y\times_U Y) $ induces a adelic divisor
\begin{align*}
\O(H)= \O(\Delta)_a - \O(\tilde{\Delta})
\end{align*}
Note that both $ G,H $ has underlying divisor $ 0 $. 
\begin{prop}
	The adelic divisors $ G,H $ are both special. 
\end{prop}

\begin{proof}
	Let $ B $ be the unit disk, and let $ p:B\to X $ be any holomorphic map with $ p(B\backslash \{0\})\subset Y $. Let $ x=p(0) $ and $ s_0=\pi(p(0)) $
	
	Then $ \pi\circ p:B\to S $ is a holomorphic map, and let $ V=(\pi\circ p)^*X $ be the pull back family on $ B $. By \cite{Zha10}, the difference of $ \omega_{X/S} $ and $ \omega_{Y/U,a} $ has the form
	\begin{align*}
	\mult_0 p^*G = g_{\mu_{ad}}(\Gamma,q;m_1,\ldots,m_r;v,v).
	\end{align*}
	Here $ m_i=p^* D_i $, and $ (\Gamma,q;m_1,\ldots,m_r) $ is the polarized metrized graph corresponding to $ V\times_S X $. We explain the point $ v\in \Gamma $ as follows. 
	
	If $ x $ is a smooth point on the fiber, then $ v $ be the vertex of the dual graph $ \Gamma $ corresponding to the unique irreducible component of $ X_{s_0} $ containing $ x $.

	If $ x $ is a node on the fiber, assume that $ D_1 $ is the boundary divisor corresponding to $ x $, then locally near $ x $, $ \pi^* D_1 $ has two components $ A,B $. Let $ a,b $ be the local intersection number of $ A,B $ with $ p(B) $. Then $ a+b=m_1 $, and if we choose coordinate $ \psi:[0,m_1]\to e $ such that the vertex $ \psi(0) $ corresponds to $ A $, the coordinate of $ v $ is just $ a $. In this way, we view $ v $ as a variable, and $ g_{\mu_{ad}}(\Gamma,q;m_1,\ldots,m_r;v,v) $ is a function of $ (r+1) $-variables. 
	
	In both cases, we see that $ G $ is special. 
	
	Similarly, if $ p:B\to T $ be any holomorphic map with $ p(B\backslash \{0\})\subset Y\times_U Y $, then by \cite{Zha10}, the adelic divisor $ H $ is characterized by 
	$$ \mult_0 p^*H = g_{\mu_{\ad}}(\Gamma,q;L_1,\ldots,L_r;v,w), $$
	where $ v,w $ are points in $ \Gamma $ given by local intersection numbers. Thus $ H $ is also special.  
\end{proof}

In this way, we may get some asymptotic behavior of Arakelov Green function and the canonical metrics, but all these results have error terms like $ o(-\sum\log |t_i|) $. José Burgos Gil, David Holmes and Robin de Jong \cite{BGHJ} have a better results. Using direct computations, they have obtained the multi-dimensional asymptotic of the canonical admissible metric on the theta line bundle on a family of complex polarized abelian varieties, see Thm. 1.1 of loc. cit. They get an O(1) error term, which is moreover continuous away from the singularities of the boundary divisor.

\bibliographystyle{plain}

\end{document}